\documentclass{amsart}%
\usepackage{amsmath,amsthm}
\usepackage{amsfonts}
\usepackage{amssymb}
\usepackage[dvips]{graphicx,color}
\usepackage{graphicx}%
\providecommand{\U}[1]{\protect\rule{.1in}{.1in}}

\begin{document}

\title{Integrable two layer point vortex motion on the half plane }
\author{M.I. Jamaloodeen\\Georgia Gwinnett College}
\date{July, 2013}
\maketitle

\begin{abstract}
In this paper we derive the equations of motion for two-layer point vortex
motion on the upper half plane. We study the invariants using symmetry,
including the Hamiltonian and show that the two vortex problem is integrable.
We characterize all two vortex motions for the cases where the vortex
strengths are both equal, $\Gamma_{1}=\Gamma_{2}=1$ and when they are opposite
$\Gamma_{1}=-\Gamma_{2}=1$. We also prove that there are no equilibria for the
two vortex problem when $\Gamma_{1}=-\Gamma_{2}=1$.\ We show that there is
only one relative equilibrium configuration when $\Gamma_{1}=\Gamma_{2}=1$ and
the vortices are in different layers. We also make observations concerning the
finite-time collapse of two vortices in the half plane. We then compare the
regimes of motion for both cases (motion on the half plane) with the case of
the two-layer vortex problem on the entire plane. We also study several
classes of streamline topologies for two vortices in different layers. We
conclude with a Hamiltonian study of integrable two-layer 3 vortex motion on
the half plane by studying integrable symmetrical configurations and provide a
rich class of new relative equilibria.

\end{abstract}

\section{\bigskip Introduction}

\noindent There is a vast literature on $N$-body vortex problems, much of
which is reference in the book by Newton \cite{Newton1}. The geostrophic
vortex models are described in many geophysics texts including \cite{Pedl}.
There has been some work done on two-layer point vortex dynamics on the entire
plane. These include the works of Young \cite{Young}, and Hogg and Stommel
\cite{HoggStommel}\cite{HoggStommell2} , and which have been primarily
numerical investigations. Some experimental work has been done by Griffiths
and Hopfinger \cite{GrifHopf}. More analytic results can be found in the work
of Gryanik \cite{Grynk} and Zabusky and McWilliams \cite{McWZab} and Flierl,
\ Polvani and Zabusky \cite{FliePolvZab}. \ Integrable two layer point vortex
motion in the plane has been extensively studied by Jamaloodeen and Newton
\cite{Jam1}\cite{JamNewt1}. More recent studies can also be found in the works
of Koshel et al \cite{KSV1}\cite{KSV2}. In these are studied the equilibrium
solutions, the vortex collapse problem and the transition to chaotic advection
through perturbations of known equilibrium solutions. The two layer vortex
problem in domains other than the entire plane has not been as extensively
studied as the one layer vortex problem on domains with boundaries. \ A good
exposition of the one layer vortex problem on domains with boundaries can be
found in work of Flucher and Gustafsson\cite{FlucGustaf}. The work in this
study can be understood as applying the techniques for one layer integrable
vortex dynamics on domains with boundaries \cite{FlucGustaf}, and the analytic
techniques for integrable two layer vortex dynamics in the unbounded plane
\cite{Jam1}\cite{JamNewt1}\cite{KSV1}\cite{KSV2}, to integrable two layer
vortex dynamics in the upper half plane.

This paper is organized as follows. We begin by deriving the equations of
motion, by first obtaining the streamfunctions for an ensemble of point
vortices. \ We also obtain the invariants through the use of symmetry and
establish the integrability of the two vortex problem in the upper half plane.
Through analysis of the Hamiltonian energy curves we characterize all 2 point
vortex motion in the upper half plane for both cases $\Gamma_{1}=$ $\Gamma
_{2}=1$ and $\Gamma_{1}=$ $-\Gamma_{2}=1$ including conclusions about 2 vortex
collapse. We then determine all equilibrium solutions for the 2 point vortex
motion, again, for both cases. We proceed by comparing the Hamiltonian energy
curves for the two layer problem in upper half plane with the one layer
problem in the upper half plane. We then present qualitative aspects of
streamline topologies for the two layer problem, for both cases $\Gamma_{1}=$
$\Gamma_{2}=1$ and $\Gamma_{1}=$ $-\Gamma_{2}=1$. \ We conclude with a study
of integrable 3 vortex motion by considering symmetrical initial
configurations, and seeking conditions to maintain the symmetry of the initial
configuration. By enforcing those conditions, which simplify by symmetry of
the configuration, we are able to obtain, numerically, rich classes of
relative equilibria in both cases $\Gamma_{1}=$ $\Gamma_{2}=1$ and $\Gamma
_{1}=$ $-\Gamma_{2}=1$.

\section{\bigskip Equations of motion and invariants}

We consider first an ensemble of two point vortices in the upper half plane
where $\Gamma_{1}$ is in bottom layer and $\Gamma_{2}$ is in top. Denoting the
respective stream functions in the corresponding layers by $\psi_{i}$
($i=1,2$) the expession for determining the streamfunctions corresponding to
unit (delta) point vortices $\delta_{i}(\varsigma)$ in each of the two layers are%

\begin{subequations}
\label{coupled}%
\begin{align}
\Delta\psi_{1}-k^{2}(\psi_{1}-\psi_{2})  &  =\delta_{1}(\zeta),\label{cup1}\\
\Delta\psi_{2}+k^{2}(\psi_{1}-\psi_{2})  &  =\delta_{2}(\zeta), \label{cup2}%
\end{align}

where the subscripts identify the layer and $1/k$ is the internal radius of
deformation. \ These equations apply here to the domain $D$, the upper half
plane with boundary $y=0$, and boundary conditions $\psi_{i}|_{\partial D}=0
$. By introducing the sum and differences $\psi=\psi_{1}+\psi_{2}$ and
$\varphi=\psi_{1}-\psi_{2}$ the equations (\ref{coupled}) uncouple%

\end{subequations}
\begin{subequations}
\label{uncup}%
\begin{align}
\Delta\psi &  =\delta_{1}+\delta_{2}=\delta(\zeta)\label{uncup1}\\
\Delta\varphi-k^{2}\varphi &  =\delta_{1}-\delta_{2}=\delta^{\prime}(\zeta)
\label{uncup2}%
\end{align}

The fundamental solution of (\ref{uncup1}) on the half plane is obtained,
through the method of images \cite{Meln}, using the Green's function,
\end{subequations}
\begin{equation}
G(z;\zeta)=G(x,y;\xi,\eta)=\ln\left\vert z-\zeta\right\vert -\ln\left\vert
z-\overline{\zeta}\right\vert \text{.} \label{Green1}%
\end{equation}
The fundamental solution of (\ref{uncup2}) on the half plane is obtained using
the Green's function
\begin{equation}
G(z;\zeta)=G(x,y;\xi,\eta)=K_{0}\left(  k\left\vert z-\zeta\right\vert
\right)  -K_{0}\left(  k\left\vert z-\overline{\zeta}\right\vert \right)
\label{Green2}%
\end{equation}
with $K_{0}\left(  x\right)  $ the modified Bessel function using the complex
variable notation,%

\[
z=x+iy\text{ and }\zeta=\xi+i\eta\text{.}%
\]

In this study we restrict the radius of deformation to be $1$ so that $k=1.$
Many regimes of motion can still be analyzed through suitable scaling coming
from initial vortex separation and/or vortex strength (circulation) assignments.

By superposition considering a vortex of strength (circulation) $\Gamma_{1}$in
layer 1 and $\Gamma_{2}$ in layer $2$ the motion of the vortices in each layer
may be computed from the streamfunctions $\psi_{1}$ and $\psi_{2}$ \ arising
from all the vortices (except the one being advected in the case of 2 or more
vortices). For example suppose there is nonzero $\Gamma_{1}$ in layer 1and
$\Gamma_{2}=0$, then stream functions induced by $\Gamma_{1}$ are,
\begin{subequations}
\begin{align}
\psi_{1}  &  =\Gamma_{1}\left(  \ln r_{1}-\ln r_{1}^{\ast}-K_{0}\left(
r_{1}\right)  +K_{0}\left(  r_{1}^{\ast}\right)  \right)  \text{,}\\
\psi_{2}  &  =\Gamma_{1}\left(  \ln r_{1}-\ln r_{1}^{\ast}+K_{0}\left(
r_{1}\right)  -K_{0}\left(  r_{1}^{\ast}\right)  \right)  \text{.}%
\end{align}

Likewise with nonzero $\Gamma_{2}$ in layer 2 and $\Gamma_{1}=0$, then stream
functions induced by $\Gamma_{2}$ are,
\end{subequations}
\begin{subequations}
\begin{align}
\psi_{1}  &  =\Gamma_{2}\left(  \ln r_{2}-\ln r_{2}^{\ast}+K_{0}\left(
r_{2}\right)  -K_{0}\left(  r_{2}^{\ast}\right)  \right) \\
\psi_{2}  &  =\Gamma_{2}\left(  \ln r_{2}-\ln r_{2}^{\ast}-K_{0}\left(
r_{2}\right)  +K_{0}\left(  r_{2}^{\ast}\right)  \right)  \text{.}%
\end{align}
Here, (see Fig. \ref{2LFig1} )
\end{subequations}
\begin{subequations}
\begin{align}
r_{i}  &  =\sqrt{\left(  x-x_{i}\right)  ^{2}+\left(  y-y_{i}\right)  ^{2}%
}\text{,}\\
r_{i}^{\ast}  &  =\sqrt{\left(  x-x_{i}\right)  ^{2}+\left(  y+y_{i}\right)
^{2}}\text{,}%
\end{align}
and use is made of
\end{subequations}
\begin{subequations}
\begin{align}
\psi_{1}  &  =1/2\left(  \psi+\varphi\right)  ,\\
\psi_{2}  &  =1/2\left(  \psi-\varphi\right)  \text{.}%
\end{align}
%

\begin{figure}
[t]
\begin{center}
\includegraphics
{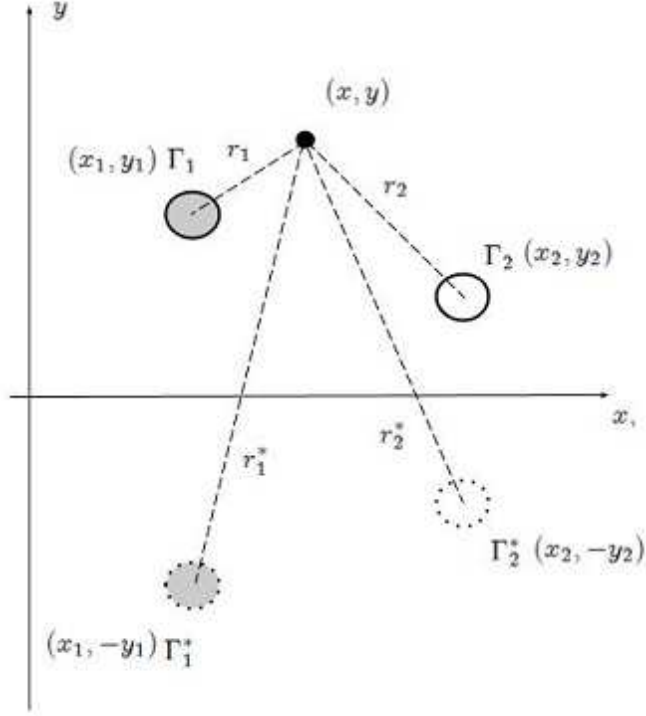}%
\caption{The basic geometry of vortices, image vortices and a tracer particle
on the half plane}%
\label{2LFig1}%
\end{center}
\end{figure}
By superposition then, with arbitrary $\varpi_{1}=\Gamma_{1}$ $\delta(r_{1})$
and $\varpi_{2}=\Gamma_{2}\delta(r_{2})$ in layers 1 and 2 respectively the
combined streamfunctions are\bigskip%
\end{subequations}
\begin{multline}
\psi_{1}=\Gamma_{1}\left(  \ln r_{1}-\ln r_{1}^{\ast}-K_{0}\left(
r_{1}\right)  +K_{0}\left(  r_{1}^{\ast}\right)  \right)  +\Gamma_{2}\left(
\ln r_{2}-\ln r_{2}^{\ast}+K_{0}\left(  r_{2}\right)  -K_{0}\left(
r_{2}^{\ast}\right)  \right)  \text{,}\label{streams_combined}\\
\psi_{2}=\Gamma_{2}\left(  \ln r_{2}-\ln r_{2}^{\ast}-K_{0}\left(
r_{2}\right)  +K_{0}\left(  r_{2}^{\ast}\right)  \right)  +\Gamma_{1}\left(
\ln r_{1}-\ln r_{1}^{\ast}+K_{0}\left(  r_{1}\right)  -K_{0}\left(
r_{1}^{\ast}\right)  \right)  \text{.}%
\end{multline}

Now the dynamics of the point vortices can be obtained by differentiation of
the streamfunctions as follows%

\begin{equation}
\dot{x}_{i}=-\left.  \frac{\partial\psi_{i}}{\partial y}\right\vert _{r=r_{i}%
},\text{ \ \ \ \ \ \ }\dot{y}_{i}=\left.  \frac{\partial\psi_{i}}{\partial
x}\right\vert _{r=r_{i}}\text{ \ \ }i=1,2\text{.} \label{dynamic1}%
\end{equation}

It is well know that the equations for point vortices are a Hamiltonian
system. \ It can be verified that the energy of the system%

\begin{equation}%
{\displaystyle\int\limits_{H}}
{\displaystyle\sum\limits_{i=1}^{2}}
\left\vert \mathbf{v}_{i}\right\vert ^{2}d\mathbf{x}\text{,}%
\end{equation}
is invariant. Integrating by parts, substituting for $\mathbf{v}_{i}$ using
(\ref{dynamic1}), $\psi_{i}$ using (\ref{streams_combined})and $\varpi_{1}$
using (\ref{uncup}) as well as the $\ 0$ boundary conditions shows that the
Hamiltonian can be simplified to,
\begin{equation}
2H(\mathbf{\varpi})=-%
{\displaystyle\int\limits_{H}}
\left(  \psi_{1}\varpi_{1}+\psi_{2}\varpi_{2}\right)  dxdy\text{,}%
\end{equation}
Using delta distributed point vortices $\varpi_{1}=\Gamma_{1}$ $\delta(r_{1}%
)$and $\varpi_{2}=\Gamma_{2}\delta(r_{2})$ and the streamfunction for the
ensemble (\ref{streams_combined}) the invariant Hamiltonian for $\Gamma_{1}$
at $\left(  x_{1},y_{1}\right)  $ and $\Gamma_{2}$ at $\left(  x_{2}%
,y_{2}\right)  $simplfies to%

\begin{multline}
H(\Gamma_{1},\Gamma_{2};\left(  x_{1},y_{1}\right)  ,\left(  x_{2}%
,y_{2}\right)  )=\label{Hamiltonian}\\
\Gamma_{1}^{2}\left[  K_{0}\left(  r_{1,1^{\ast}}\right)  -\ln r_{1,1^{\ast}%
}\right]  +2\Gamma_{1}\Gamma_{2}\left[  \ln r_{1,2}-\ln r_{1,2^{\ast}}%
+K_{0}\left(  r_{1,2}\right)  -K_{0}\left(  r_{1,2^{\ast}}\right)  \right]
+\Gamma_{2}^{2}\left[  K_{0}\left(  r_{2,2^{\ast}}\right)  -\ln r_{2,2^{\ast}%
}\right]  \text{,}%
\end{multline}
again with reference to Fig.(\ref{2LFig1}). Image vortices are denoted with a
$\ast$. Note also from the symmetry of the geometry of Fig. (\ref{2LFig1})
that $r_{1,2^{\ast}}=r_{2,1^{\ast}}$ and $r_{1,2}=r_{1^{\ast},2^{\ast}}$.

When $\Gamma_{2}=0$, $\ $corresponding to a single vortex, the Hamiltonian,
$H$ simplifies to $H=\Gamma_{1}^{2}\left[  K_{0}\left(  r_{1,1^{\ast}}\right)
-\ln r_{1,1^{\ast}}\right]  =C$ and by the monotonicity of $K_{0}\left(
r_{1,1^{\ast}}\right)  -\ln r_{1,1^{\ast}}$ we conclude that $r_{1,1^{\ast}%
}=constant$ corresponding to the vortex $\Gamma_{1}$ translating parallel to
the $x$-axis. \ This is the general solution then for the case of a single
vortex. Notice also in (\ref{Hamiltonian}) that the Hamiltonian, \ $H$, is
invariant with respect to arbitrary displacements of both $\Gamma_{1}$ at
$\left(  x_{1},y_{1}\right)  $ and $\Gamma_{2}$ at $\left(  x_{2}%
,y_{2}\right)  $ by $\delta$ in the $x$-direction. This implies by Noether's
theorem the invariance of%

\begin{equation}
\Gamma_{1}y_{1}+\Gamma_{2}y_{2}=c\text{.} \label{y_momentum}%
\end{equation}

We provide an explicit proof of (\ref{y_momentum}) in the appendix. The
invariance of the Hamiltonian, \ $H$, (\ref{Hamiltonian}) and the momentum in
the $y$-direction (\ref{y_momentum}) imply that the 2 layer, 2 point vortex
system in the half-plane is integrable. An easy way to see this is to notice
that the Hamiltonian \ depends only on $r_{1,1^{\ast}}=2y_{1}$ and
$r_{2,2^{\ast}}=2y_{2}$ so that by (\ref{y_momentum}) it depends only on
either $y_{1}$ or $y_{2}$. The other term in the Hamiltonian is $r_{1,2}%
=\sqrt{\left(  x_{1}-x_{2}\right)  ^{2}+\left(  y_{1}-y_{2}\right)  ^{2}}$,
from which we conclude that the Hamiltonian is a function of the two
variables, $y_{1}$ (or $y_{2}$) and $\left\vert x_{1}-x_{2}\right\vert $.

\section{\bigskip Characterizing 2 point vortex motion}

\noindent We consider the cases $\Gamma_{1}=$ $\Gamma_{2}=1$ and $\Gamma_{1}=$
$-\Gamma_{2}=1$ separately.

\subsection{The case $\Gamma_{1}=$ $-\Gamma_{2}=1$}

In this case we use the invariant $\Gamma_{1}y_{1}+\Gamma_{2}y_{2}=c$ or
$y_{1}-y_{2}=\alpha$ with parameter $\alpha$. With reference to Fig.
(\ref{2LFig1}) \ we designate $y_{1}=y$ so that $y_{2}=y-\alpha$ (so that
$y_{1}^{\ast}=-y$ and $y_{2}^{\ast}=\alpha-y$ and denote $x=\left\vert
x_{1}-x_{2}\right\vert $ with the Hamiltonion (\ref{Hamiltonian}) simplifying
to%
\begin{multline}
H(x,y)=K_{0}\left(  2y\right)  +K_{0}\left(  2\left(  y-\alpha\right)
\right)  -\ln\left(  2y\right)  \left(  2\left(  y-\alpha\right)  \right)
\label{phase1}\\
+\ln\frac{x^{2}+\left(  2y-\alpha\right)  ^{2}}{x^{2}+\alpha^{2}}%
+2K_{0}\left(  \sqrt{x^{2}+\left(  2y-\alpha\right)  ^{2}}\right)
-2K_{0}\left(  \sqrt{x^{2}+\alpha^{2}}\right)  .
\end{multline}

\begin{figure}
[ptbh]
\begin{center}
\includegraphics[
natheight=8.260700in,
natwidth=5.750100in,
height=6.2241in,
width=4.3405in
]%
{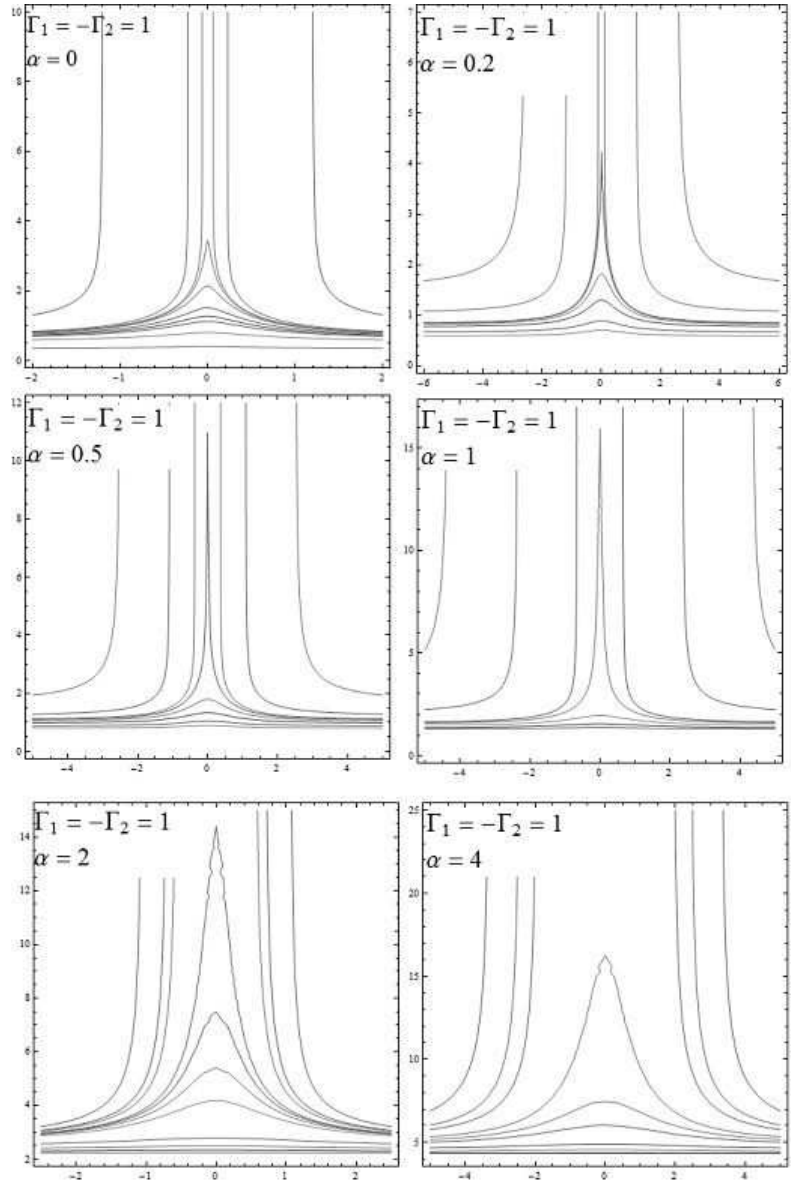}%
\caption{Hamiltonian level curves for $\Gamma_{1}=-\Gamma_{2}=1$, and
parameter values $\alpha=0,0.2,0.5,1,2,4$.}%
\label{2LFig2}%
\end{center}
\end{figure}

The Hamiltonian level curves are shown for various parameter values $\alpha
$\ in Fig. (\ref{2LFig2}). There are no periodic solutions, and, in
particular, no equilibria. We shall rigorously justifiy this in the sequel.
There are two types of motion. The first corresponds to $\left\vert
x_{1}-x_{2}\right\vert \rightarrow\infty$ and $y_{1}=y_{2}$ approaching a
nonzero value, or vice versa $y_{1}=y_{2}\rightarrow\infty$ with $x=\left\vert
x_{1}-x_{2}\right\vert $ approaching a nonzero value. \ Which of these will
depend on which side of the phase plane one begins; either on the front side
or the back side of the phase plane. This is similar to the only type of
motion seen for the one layer problem. See Fig (\ref{2LFig4}) (a), for the
corresponding case\ $y_{1}-y_{2}=\alpha=0$ for the one layer problem in the
half plane. There is a second type of motion for the two layer problem seen,
for all values of $\alpha$ including $\alpha=0$. These are phase curves that
cross the $x=0$ $\ $\ or $\left\vert x_{1}-x_{2}\right\vert =0$ line. In the
case that $\alpha=0$ this crossing corresponds to $\alpha=0=y_{1}-y_{2\text{
}}$ and $\left\vert x_{1}-x_{2}\right\vert =0$, which corresponds to vortex
collapsing configurations. We are still investigating the nature of these
collapsing configurations, as to whether they are finite time or infinite time
vortex collapse solutions. Preliminary numerical results suggest that these
are infinite time collapsing configurations. Notice that these collapsing
configurations are not admitted for $\Gamma_{1}=$ $-\Gamma_{2}=1$, in the one
layer case in the upper half plane. \ See Fig (\ref{2LFig4}) (a), for the
corresponding case\ $y_{1}-y_{2}=\alpha=0$ for the one layer problem in the
half plane and Fig (\ref{2LFig4}) (b) for the case $y_{1}-y_{2}=\alpha=1$.
\ The second type of motion is seen in the one layer case, provide $\alpha>0$
such as $y_{1}-y_{2}=\alpha=1$, however while $\left\vert x_{1}-x_{2}%
\right\vert \rightarrow0$, the $y$ values do not since $y_{1}-y_{2}=\alpha=1$.%

\bigskip

\subsection{The case $\Gamma_{1}=$ $\Gamma_{2}=1$}

In this case we use the invariant $\Gamma_{1}y_{1}+\Gamma_{2}y_{2}=c$ or
$y_{1}+y_{2}=\alpha$ with parameter $\alpha$. With reference to Fig.
(\ref{2LFig1}) \ we designate $y_{1}=y$ so that $y_{2}=\alpha-y$ and denote
$x=\left\vert x_{1}-x_{2}\right\vert $ with the Hamiltonian (\ref{Hamiltonian}%
) simplifying to%
\begin{multline}
H(x,y)=K_{0}\left(  2y\right)  +K_{0}\left(  2\left(  \alpha-y\right)
\right)  -\ln\left(  2y\right)  \left(  2\left(  \alpha-y\right)  \right)
\label{phase2}\\
+\ln\frac{x^{2}+\left(  2y-\alpha\right)  ^{2}}{x^{2}+\alpha^{2}}%
+2K_{0}\left(  \sqrt{x^{2}+\left(  2y-\alpha\right)  ^{2}}\right)
-2K_{0}\left(  \sqrt{x^{2}+\alpha^{2}}\right)  .
\end{multline}

\begin{figure}
[ptb]
\begin{center}
\includegraphics[
natheight=8.291800in,
natwidth=5.666200in,
height=6.2465in,
width=4.2782in
]%
{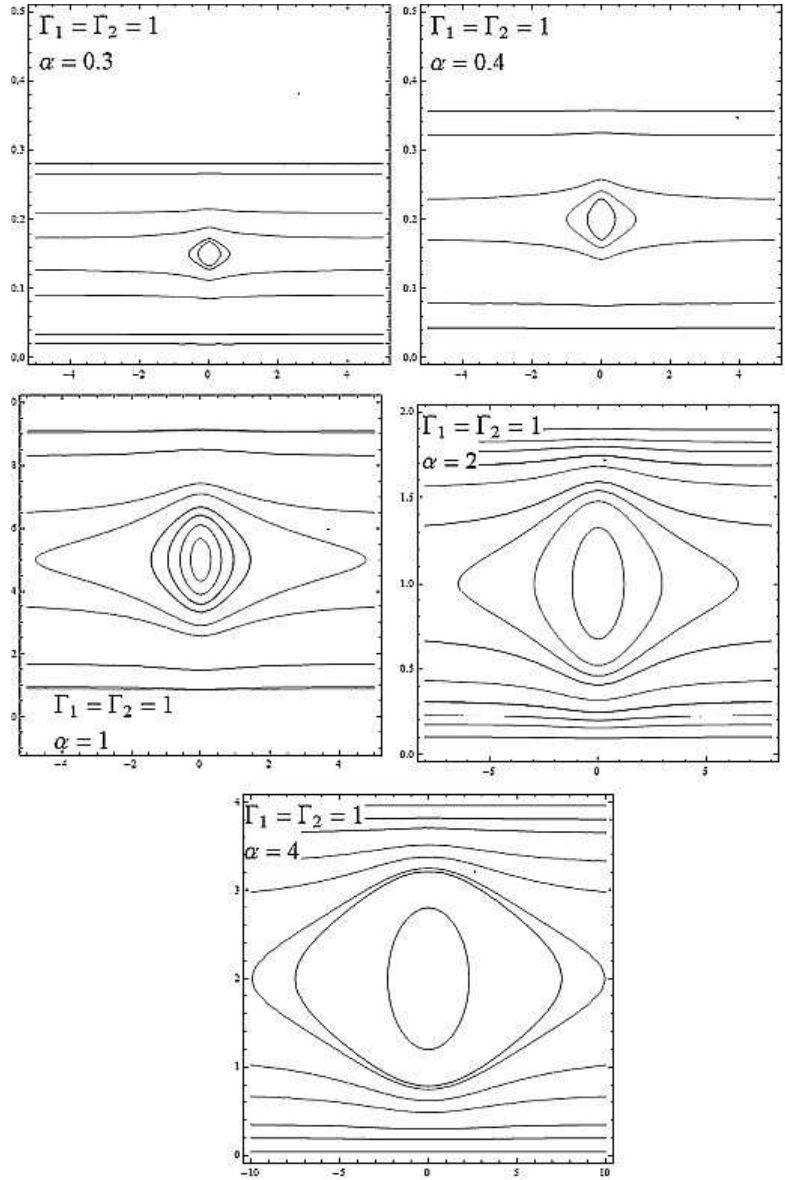}%
\caption{Hamiltonian level curves for $\Gamma_{1}=\Gamma_{2}=1$, and parameter
values $\alpha=0,0.3,0.4,1,2,4$.}%
\label{2LFig3}%
\end{center}
\end{figure}

\bigskip The Hamiltonian level curves are shown for various parameter values
$\alpha$\ in Fig. (\ref{2LFig3}). The motion in this case is very similar in
many ways to the one layer problem with $\Gamma_{1}=$ $\Gamma_{2}=1$,(see Fig
(\ref{2LFig4}) (c)) with one minor difference. \ The similarities include the
two types of motions. Unlike the case $\Gamma_{1}=$ $\Gamma_{2}=-1$, the case
$\Gamma_{1}=$ $\Gamma_{2}=1$ admits closed periodic solutions. Also in the two
layer case the center of these periodic curves corresponds to a fixed
equilibrium. As previously mentioned the two layer problem in the upper half
plane for $\Gamma_{1}=$ $\Gamma_{2}=-1$ does not admit any equilibrium
solutions. What is interesting about these fixed equilibria is that they are
centered at $\left(  x,y\right)  =(0,\frac{\alpha}{2})$ corresponding to
$x_{1}=x_{2\text{ }}$ and $y_{1}=y_{2}=\frac{\alpha}{2}$. The coordinates of
both vortices are the same meaning they are stacked one on top of the other
for the fixed equilibrium configuration. \ We will study, more carefully, in
what follows this equilibrium solution and rigorously show that there are no
other equilibrium solutions for the case $\Gamma_{1}=$ $\Gamma_{2}=1$ in the
upper half plane.

Clearly for the one layer problem we cannot have both $x_{1}=x_{2\text{ }}$and
$y_{1}=y_{2}$ since there is freedom to stack the vortices, and so this would
correspond to a collapsing configuration in which the vortices where initially
located such that $x_{1}=x_{2\text{ }}$and $y_{1}=y_{2}$ which is not
feasible. This is the one major difference between the two layer and one layer
problems in the upper half plane for the case $\Gamma_{1}=$ $\Gamma_{2}=1$.

The second type of motion observed in Fig. (\ref{2LFig3}), is a non periodic
regime of motion in which either $\left\vert x_{1}-x_{2}\right\vert
\rightarrow0$ and $y_{1}$ approaches a nonzero number or where $\left\vert
x_{1}-x_{2}\right\vert \rightarrow\infty$ and $y_{1}$ approaches a nonzero
number. Which of these will again depend on which side of the phase plane one
begins; either on the front side or the back side of the phase plane.

\subsection{Equilibrium solutions}

\begin{figure}
[t!]
\begin{center}
\includegraphics[
natheight=5.521000in,
natwidth=4.791900in,
height=3.8925in,
width=3.3823in
]%
{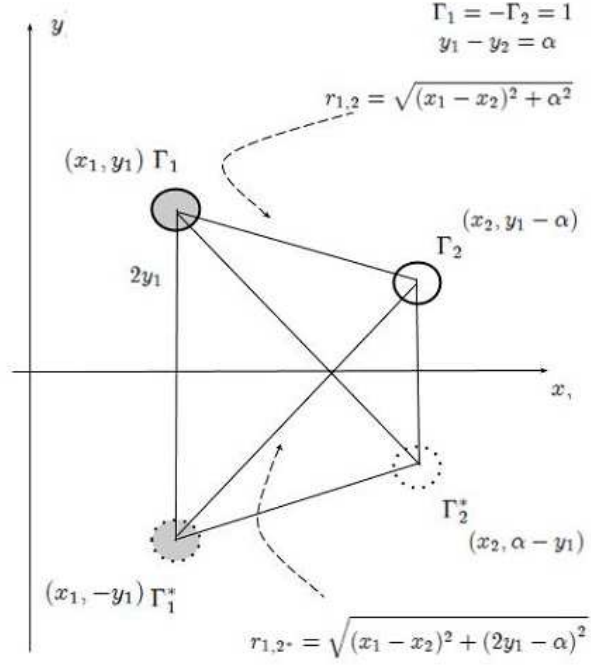}%
\caption{Two vortices in different layers $\Gamma_{1}=$ -$\Gamma_{2}=1$. In
this case the invariant $\Gamma_{1}y_{1}+\Gamma_{2}y_{2}=\alpha$, becomes
$y_{1}-y_{2}=\alpha$. $\Gamma_{1}$ in layer 1 at $\left(  x_{1,}y_{1}\right)
$, $\Gamma_{2}$ in layer 2 at $\left(  x_{2,}y_{1}-\alpha\right)  $ and
$r_{1,2}=\sqrt{\left(  x_{1}-x_{2}\right)  ^{2}+\alpha^{2}}$ and
$r_{1,2^{\ast}}=\sqrt{\left(  x_{1}-x_{2}\right)  ^{2}+\left(  2y_{1}%
-\alpha\right)  ^{2}}$.}%
\label{2L2v_g1_negative_g2}%
\end{center}
\end{figure}

In this section we consider equilibrium solutions for the 2 vortex problem. We
begin by showing that there are no relative equilibrium solutions when
$\Gamma_{1}=$ $-\Gamma_{2}=1$ with vortex $1$, $\Gamma_{1}$ at $\left(
x_{1},y_{1}\right)  $ in layer $1$ and vortex $2$, $\Gamma_{2}$ at $\left(
x_{2},y_{2}\right)  =\left(  x_{2},y_{1}-\alpha\right)  $ in layer $2$.\ In
this case the invariant $\Gamma_{1}y_{1}+\Gamma_{2}y_{2}=\alpha$, becomes
$y_{1}-y_{2}=\alpha$. We consider the distance (see Fig.
(\ref{2L2v_g1_negative_g2})) between the two vortices $\Gamma_{1}$ and
$\Gamma_{2}$, $r_{1,2}(t)=\sqrt{\left(  x_{1}-x_{2}\right)  ^{2}+\alpha^{2}}$
and show that it is never constant. Consider,%

\begin{equation}
\frac{dr_{1,2}^{2}}{dt}=\left(  x_{1}-x_{2}\right)  \left(  \dot{x}_{1}%
-\dot{x}_{2}\right)  +\left(  y_{1}-y_{2}\right)  \left(  \dot{y}_{1}-\dot
{y}_{2}\right)  =\left(  x_{1}-x_{2}\right)  \left(  \dot{x}_{1}-\dot{x}%
_{2}\right)  +\left(  y_{1}-y_{2}\right)  (0)=\left(  x_{1}-x_{2}\right)
\left(  \dot{x}_{1}-\dot{x}_{2}\right)  \text{,}%
\end{equation}

where the last follows by the invariance of $y_{1}-y_{2}=\alpha$ or $\dot
{y}_{1}-\dot{y}_{2}=0$. Using the equations (\ref{dynamicx1},\ref{dynamicx2})
and $\Gamma_{1}=-$ $\Gamma_{2}=1$, we obtain,%

%

\begin{equation}
\left(  \dot{x}_{1}-\dot{x}_{2}\right)  =K_{1}(2y_{1})+K_{1}(2y_{2}%
)+\frac{2(y_{1}+y_{2})K_{1}(r_{1,2^{\ast}})}{r_{1,2^{\ast}}}+\frac{1}{2y_{1}%
}+\frac{1}{2y_{2}}-\frac{2(y_{1}+y_{2})}{r_{1,2^{\ast}}^{2}}>0\text{,}%
\end{equation}
where the inequality $\left(  \dot{x}_{1}-\dot{x}_{2}\right)  >0$ follows from
the positivity $y_{1}>0,y_{2}>0,K_{1}(x)>0$ and the fact that
\begin{equation}
\frac{1}{2y_{1}}+\frac{1}{2y_{2}}-\frac{2(y_{1}+y_{2})}{r_{1,2^{\ast}}^{2}%
}>0\text{,}%
\end{equation}

which can easily shown using the geometry of the upper half plane.

\begin{figure}
[h]
\begin{center}
\includegraphics[
natheight=5.521000in,
natwidth=4.791900in,
height=4.1693in,
width=3.6218in
]%
{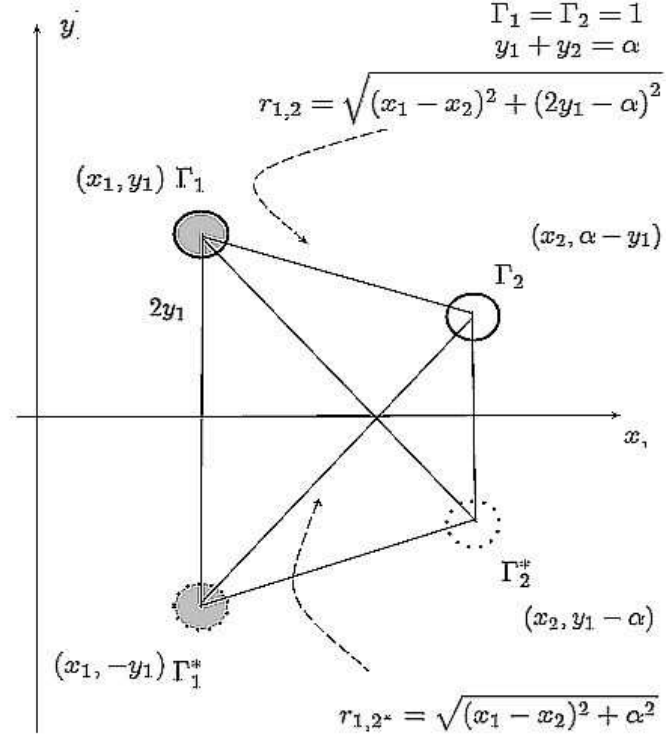}%
\caption{Two vortices in different layers $\Gamma_{1}=$ $\Gamma_{2}=1$. In
this case the invariant $\Gamma_{1}y_{1}+\Gamma_{2}y_{2}=\alpha$, becomes
$y_{1}+y_{2}=\alpha$. $\Gamma_{1}$ in layer 1 at $\left(  x_{1,}y_{1}\right)
$, $\Gamma_{2}$ in layer 2 at $\left(  x_{2,}y_{1}-\alpha\right)  $ and
$r_{1,2}=\sqrt{\left(  x_{1}-x_{2}\right)  ^{2}+\left(  2y_{1}-\alpha\right)
^{2}}$ and $r_{1,2^{\ast}}=\sqrt{\left(  x_{1}-x_{2}\right)  ^{2}+\alpha^{2}}%
$.}%
\label{2L2v_g1_positive_g2}%
\end{center}
\end{figure}

Next we show that there is only one relative equilibrium solutions when
$\Gamma_{1}=$ $\Gamma_{2}=1$ with vortex $1$, $\Gamma_{1}$ at $\left(
x_{1},y_{1}\right)  $ in layer $1$ and vortex $2$, $\Gamma_{2}$ at $\left(
x_{2},y_{2}\right)  =\left(  x_{2},\alpha-y_{1}\right)  $ in layer $2$.\ This
corresponds to the case $\left(  x_{1,}y_{1}\right)  =\left(  x_{2,}%
y_{2}\right)  $ corresponding to two vortices lying exactly one on top of the
other. \ Note this equilibrium position is not feasible in the one-layer case.
These relative equilibria are clear in Fig. (\ref{2LFig3}) in which
$x=x_{1}-x_{2}=0$, and $y_{1}=y_{2}=\frac{\alpha}{2}$ and correspond to the
center of the periodic orbits. In this case the invariant $\Gamma_{1}%
y_{1}+\Gamma_{2}y_{2}=\alpha$, becomes $y_{1}+y_{2}=\alpha$. We consider the
distance (see Fig.(\ref{2L2v_g1_positive_g2})) between the two vortices
$\Gamma_{1}$ and $\Gamma_{2}$, $r_{1,2}(t)=\sqrt{\left(  x_{1}-x_{2}\right)
^{2}+\left(  2y_{1}-\alpha\right)  ^{2}}$ and show that it is never constant.%

Consider,%

\begin{equation}
\left(  \dot{x}_{1}-\dot{x}_{2}\right)  =K_{1}(2y_{1})-K_{1}(2y_{2}%
)+\frac{2(y_{1}-y_{2})K_{1}(r_{1,2})}{r_{1,2}}+\frac{1}{2y_{1}}-\frac
{1}{2y_{2}}-\frac{2(y_{1}-y_{2})}{r_{1,2}^{2}}\text{,}
\label{x1_minus_x2_both_positive}%
\end{equation}

\begin{equation}
\dot{y}_{1}=\left(  x_{1}-x_{2}\right)  \left(  \frac{K_{1}(r_{1,2^{\ast}}%
)}{r_{1,2^{\ast}}}-\frac{K_{1}(r_{1,2})}{r_{1,2}}+\frac{1}{r_{1,2}^{2}}%
-\frac{1}{r_{1,2^{\ast}}^{2}}\right)  \text{,}%
\end{equation}

\begin{multline}
\frac{dr_{1,2}^{2}}{dt}=\left(  x_{1}-x_{2}\right)  \left(  \dot{x}_{1}%
-\dot{x}_{2}\right)  +\left(  y_{1}-y_{2}\right)  \left(  \dot{y}_{1}-\dot
{y}_{2}\right)  =\label{relativeequation_both_positive}\\
\left(  x_{1}-x_{2}\right)  \left(  K_{1}(2y_{1})-K_{1}(2y_{2})+\frac
{2(y_{1}-y_{2})K_{1}(r_{1,2^{\ast}})}{r_{1,2^{\ast}}}+\frac{1}{2y_{1}}%
-\frac{1}{2y_{2}}-\frac{2(y_{1}-y_{2})}{r_{1,2^{\ast}}^{2}}\right)  \text{.}%
\end{multline}

By considering the phase curves $\Gamma_{1}=$ $\Gamma_{2}=1$ shown in Fig.
(\ref{2LFig3}) we see that the curves all pass through $x=(x_{1}-x_{2})=0$.
\ We show that in this case $\left(  \dot{x}_{1}-\dot{x}_{2}\right)  \neq0$ so
that in Eq. (\ref{relativeequation_both_positive}) for a relative equilibrium%

\begin{equation}
K_{1}(2y_{1})-K_{1}(2y_{2})+\frac{2(y_{1}-y_{2})K_{1}(r_{1,2^{\ast}}%
)}{r_{1,2^{\ast}}}+\frac{1}{2y_{1}}-\frac{1}{2y_{2}}-\frac{2(y_{1}-y_{2}%
)}{r_{1,2^{\ast}}^{2}}=0\text{.} \label{relativeequation_both_positive_ypart}%
\end{equation}

Now when $x=x_{1}-x_{2}=0$, (\ref{x1_minus_x2_both_positive}) becomes%

\begin{equation}
\left(  \dot{x}_{1}-\dot{x}_{2}\right)  =K_{1}(2y_{1})-K_{1}(2y_{2}%
)+\frac{2(y_{1}-y_{2})K_{1}(\left\vert y_{1}-y_{2}\right\vert )}{\left\vert
y_{1}-y_{2}\right\vert }+\frac{1}{2y_{1}}-\frac{1}{2y_{2}}-\frac{2(y_{1}%
-y_{2})}{\left\vert y_{1}-y_{2}\right\vert ^{2}}\text{,}%
\end{equation}

\bigskip

which can be shown to be nonzero except when $y_{1}-y_{2}=0$ which would be
simultaneous with $x_{1}-x_{2}=0$ So in Eq.
(\ref{relativeequation_both_positive}) we require
(\ref{relativeequation_both_positive_ypart}). In this case when $x_{1}%
-x_{2}=0$, (\ref{relativeequation_both_positive_ypart}) simplifies to

\bigskip%
\begin{equation}
K_{1}(2y_{1})-K_{1}(2y_{2})+\frac{2(y_{1}-y_{2})K_{1}(y_{1}+y_{2}%
)}{r_{1,2^{\ast}}}+\frac{1}{2y_{1}}-\frac{1}{2y_{2}}-\frac{2(y_{1}-y_{2}%
)}{(y_{1}+y_{2})^{2}}=0\text{,}%
\end{equation}

which again can be shown to be nonzero except when $y_{1}-y_{2}=0$ which would
be simultaneous with $x_{1}-x_{2}=0$. The case $\Gamma_{1}=$ $\Gamma_{2}=1$
admits then the relative equilibria $\left(  x_{1,}\frac{\alpha}{2}\right)
=\left(  x_{2,}\frac{\alpha}{2}\right)  $ where the vortices lie exactly one
on top of the other. \ We can see that the phase curves (\ref{phase2}) near
this equilibrium position are closed periodic orbits so that the phase curves
near to $\left(  x_{1,}\frac{\alpha}{2}\right)  =\left(  x_{2,}\frac{\alpha
}{2}\right)  $ are almost relative equilibria, in the sense that the level
curves (\ref{phase2}) close to $\left(  x_{1,}\frac{\alpha}{2}\right)
=\left(  x_{2,}\frac{\alpha}{2}\right)  $ are asymptotically elliptical
$\left(  x_{1}-x_{2}\right)  ^{2}+(y_{1}-y_{2})^{2}=c^{2}$ or $\left(
x_{1}-x_{2}\right)  ^{2}+(y_{1}-(\alpha-y_{1}))^{2}=c^{2}$ or $\left(
x_{1}-x_{2}\right)  ^{2}+(2y_{1}-\alpha))^{2}=c^{2}$ corresponding to what
would be a true relative equilibrium.\bigskip\ \ This is a novel relative
equilibrium solution, keeping in mind that the one-layer two vortex problem on
the upper half plan admits no relative equilibrium solutions.

\subsection{Comparing with the one layer two vortex problem in the upper half
plane}

We summarize the one layer two vortex regimes of motion in the upper half
plane for completeness and to highlight the similarities and differences we
mentioned earlier.

In this case the Green's functions using the method of images corresponds only
to solving in one layer (\ref{uncup1}) and not (\ref{uncup2}) and is given by
(\ref{Green1}). \ The general Hamiltonian becomes,%
\begin{equation}
H(\Gamma_{1},\Gamma_{2};\left(  x_{1},y_{1}\right)  ,\left(  x_{2}%
,y_{2}\right)  )=-\Gamma_{1}^{2}\ln r_{1,1^{\ast}}+2\Gamma_{1}\Gamma_{2}%
\ln\frac{r_{1,2}}{r_{1,2^{\ast}}}-\Gamma_{2}^{2}\ln r_{2,2^{\ast}}\text{.}%
\end{equation}

Again, use is made of the invariant $\Gamma_{1}y_{1}+$ $\Gamma_{2}y_{2}%
=\alpha$. We consider three representative cases:

\begin{enumerate}
\item $\Gamma_{1}=-\Gamma_{2}=1$ and $y_{1}-$ $y_{2}=\alpha=0.$ The
Hamiltonian , using $x=\left\vert x_{1}-x_{2}\right\vert $ and $y_{1}=y_{2}=y$,%

\begin{equation}
\frac{H(x,y)}{2}=\ln\frac{\sqrt{x^{2}+\left(  2y\right)  ^{2}}}{2yx}\text{.}%
\end{equation}

\item $\Gamma_{1}=-\Gamma_{2}=1$ and $y_{1}-$ $y_{2}=\alpha=1$, so
$y_{2}=y_{1}-\alpha=y-1$, with Hamiltonian%

\begin{equation}
H(x,y)=-\ln(4y(y-1))+2\ln\frac{\sqrt{x^{2}+\left(  2y-1\right)  ^{2}}}%
{\sqrt{x^{2}+1^{2}}}\text{.}%
\end{equation}

\item $\Gamma_{1}=\Gamma_{2}=1$ and $y_{1}+y_{2}=\alpha=1,$so $y_{2}%
=\alpha-y_{1}=1-y$, with Hamiltonian%

\begin{equation}
H(x,y)=-\ln(4y(1-y))+2\ln\frac{\sqrt{x^{2}+\left(  2y-1\right)  ^{2}}}%
{\sqrt{x^{2}+1^{2}}}\text{.}%
\end{equation}

The phase curves are shown in Fig.(\ref{2LFig4}) (a), (b) and (c) respectively.%

\begin{figure}
[h]
\begin{center}
\includegraphics[
natheight=5.500200in,
natwidth=5.719000in,
height=4.4278in,
width=4.6017in
]%
{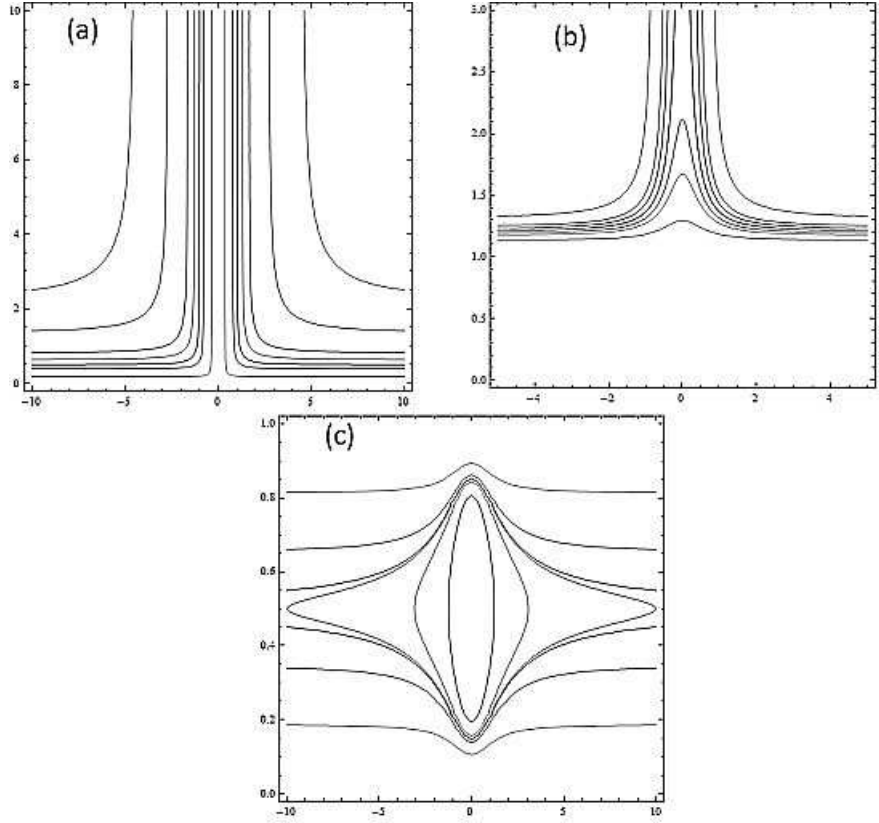}%
\caption{Hamiltonian phase plots for the one layer two vortex problem
(a)$\Gamma=-\Gamma_{2}=1$ and $y_{1}-$ $y_{2}=\alpha=0.$ (b)$\Gamma
=-\Gamma_{2}=1$ and $y_{1}-$ $y_{2}=\alpha=1$, and (c) $\Gamma=\Gamma_{2}=1$
and $y_{1}+y_{2}=\alpha=1$.}%
\label{2LFig4}%
\end{center}
\end{figure}

Notice, in particular, that when $\Gamma_{2}=0$ $\ $corresponding to a single
vortex, the Hamiltonian, $H$ simplifies to $H=\Gamma_{1}^{2}\left[
K_{0}\left(  r_{1,1^{\ast}}\right)  -\ln r_{1,1^{\ast}}\right]  =C$ and by the
monotonicity of $K_{0}\left(  r_{1,1^{\ast}}\right)  -\ln r_{1,1^{\ast}}$ we
conclude that $r_{1,1^{\ast}}=constant$ corresponding to the vortex
$\Gamma_{1}$ translating parallel to the $x$-axis. \ This is the general
solution then for the case of a single vortex. Notice also in
(\ref{Hamiltonian}) that the Hamiltonian, \ $H$, is invariant with respect to
arbitrary displacements of both $\Gamma_{1}$ at $\left(  x_{1},y_{1}\right)  $
and $\Gamma_{2}$ at $\left(  x_{2},y_{2}\right)  $ by $\delta$ in the
$x$-direction. This implies by Noether's theorem the invariance of
\end{enumerate}

\begin{equation}
\Gamma_{1}y_{1}+\Gamma_{2}y_{2}=c\text{.}%
\end{equation}

.

\section{\bigskip Streamline topologies for the two vortex problem.}

\noindent We consider the cases $\Gamma_{1}=$ $\Gamma_{2}=1$ and $\Gamma_{1}=$
$-\Gamma_{2}=1$ separately.

\subsection{The case $\Gamma_{1}=$ $\Gamma_{2}=1$}

\begin{figure}
[h]
\begin{center}
\includegraphics
{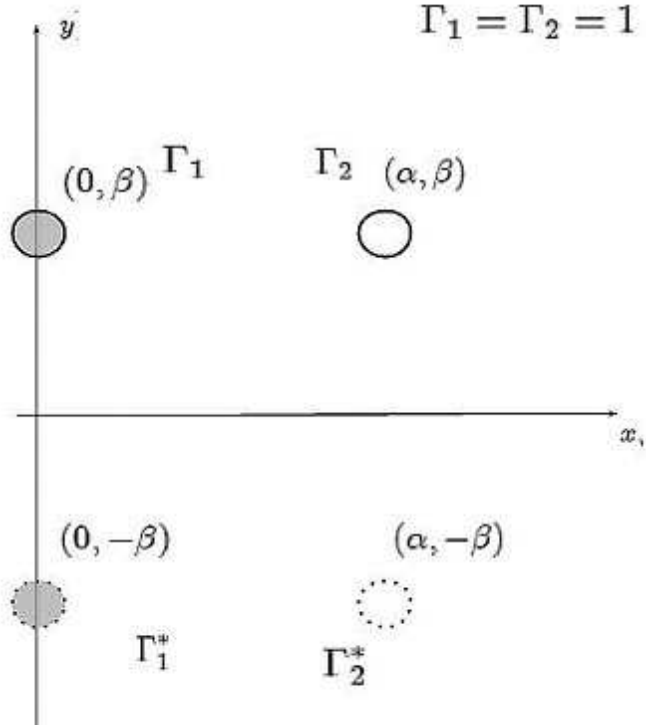}%
\caption{The configuration for streamline topologies with parameters $\alpha$,
$\beta$. The streamline is $\psi_{1}=\ln\sqrt{x^{2}+(y-\beta)^{2}}%
-K_{0}\left(  \sqrt{x^{2}+(y-\beta)^{2}}\right)  +K_{0}\left(  \sqrt
{x^{2}+(y+\beta)^{2}}\right)  -\ln\sqrt{x^{2}+(y+\beta)^{2}}+\ln
\sqrt{(x-\alpha)^{2}+(y-\beta)^{2}}+K_{0}\left(  \sqrt{(x-\alpha)^{2}%
+(y-\beta)^{2}}\right)  -K_{0}\left(  \sqrt{(x-\alpha)^{2}+(y+\beta)^{2}%
}\right)  -\ln\left(  \sqrt{(x-\alpha)^{2}+(y+\beta)^{2}}\right)  $}%
\label{Stream_configuration1}%
\end{center}
\end{figure}

We consider two basic configurations. One in which $\Gamma_{1}=$ $\Gamma
_{2}=1$, with $\Gamma_{1}\ $\ in layer $1$, and $\Gamma_{2}\ $\ in layer $2 $,
with the same $y$-coordinate ($\beta$) as shown in Fig
(\ref{Stream_configuration1}). We consider various values of the parameters
$\alpha$ and $\beta$ with reference to Rossby radius of deformation $k=1$.%

\begin{figure}
[h]
\begin{center}
\includegraphics[
natheight=5.124900in,
natwidth=5.083400in,
height=5.1802in,
width=5.1387in
]%
{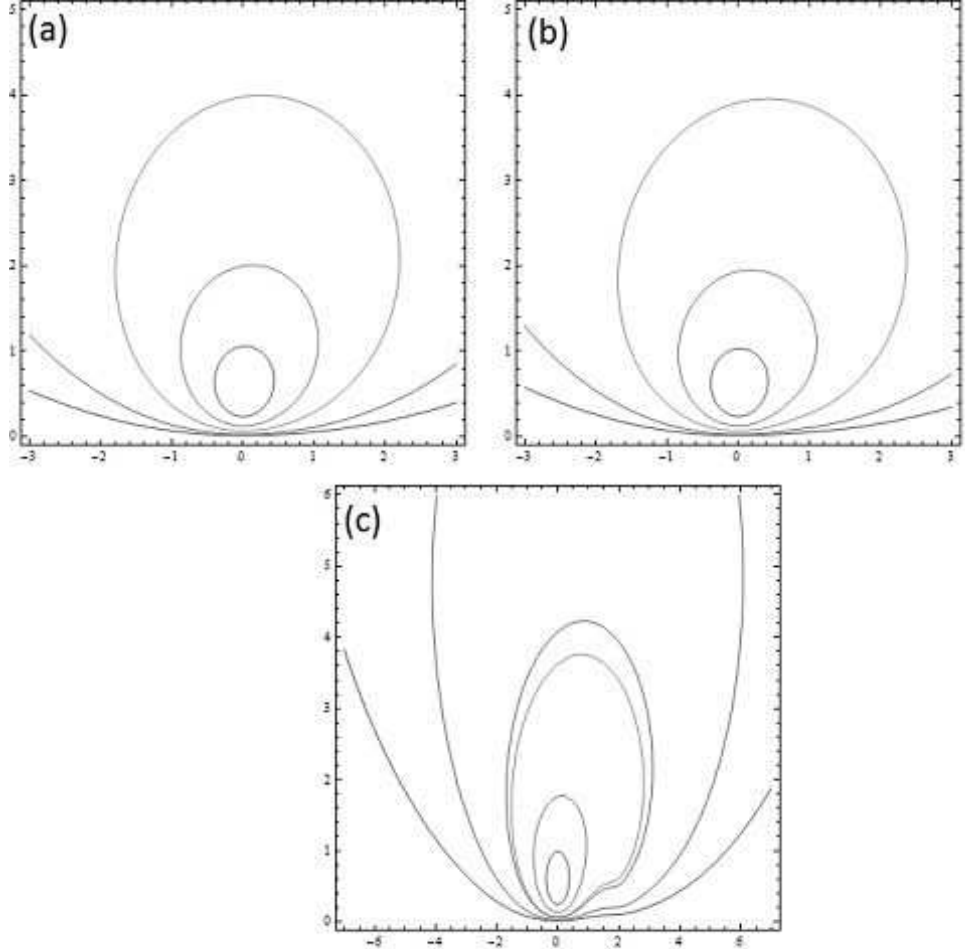}%
\caption{Streamlines for horizontal configurations. $\Gamma_{1}=\Gamma_{2}$ in
layers $1$ and $2$ respectively. $\beta=0.5$. (a) $\alpha=0.6$, (b) $\alpha
=1$, and (c) $\alpha=2$.}%
\label{Streamline1}%
\end{center}
\end{figure}

\begin{figure}
[h!]
\begin{center}
\includegraphics[
natheight=5.666200in,
natwidth=5.572800in,
height=5.6948in,
width=5.6005in
]%
{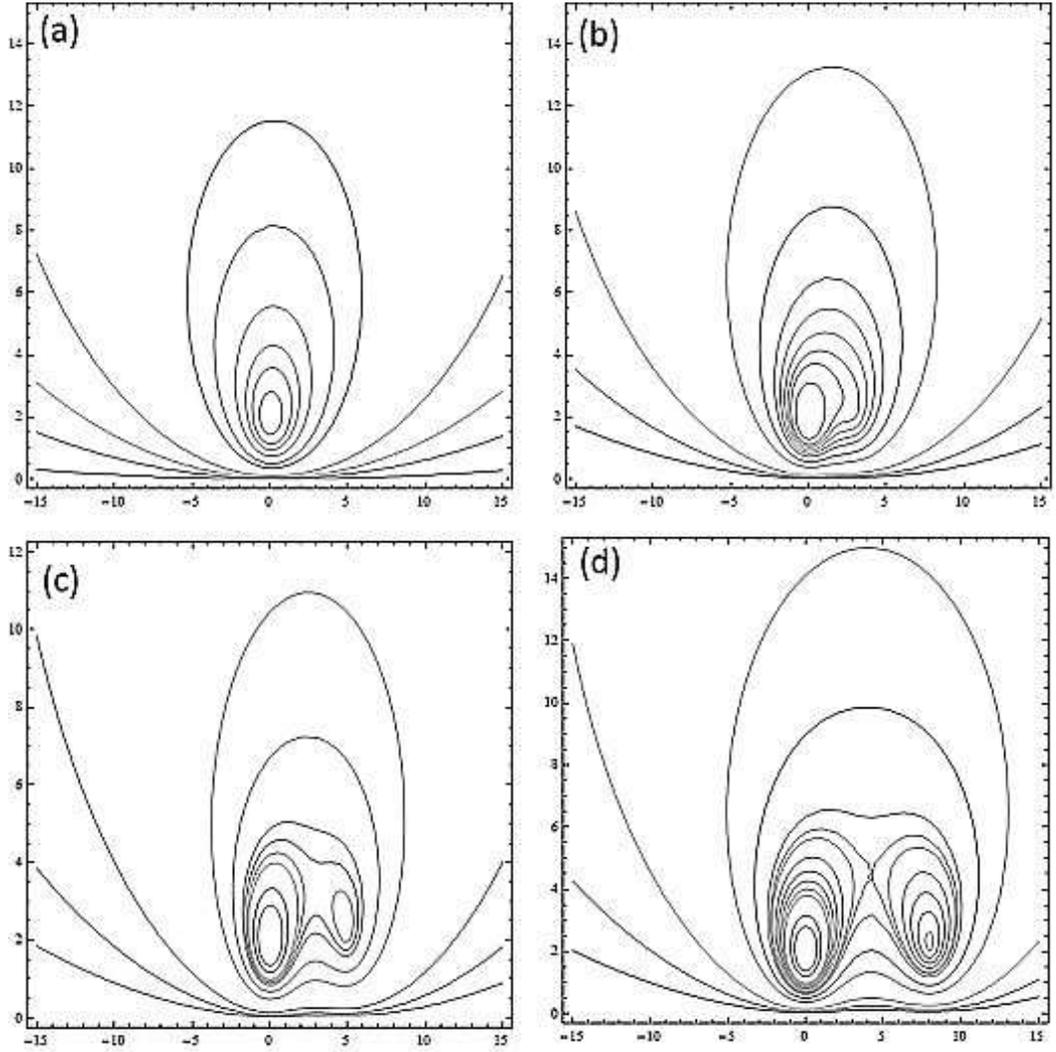}%
\caption{Streamlines for horizontal configurations. $\Gamma_{1}=\Gamma_{2}$ in
layers $1$ and $2$ respectively. $\beta=2$. (a) $\alpha=0.6$, (b) $\alpha=3$,
(c) $\alpha=5$ and (d) $\alpha=8$}%
\label{Streamline2}%
\end{center}
\end{figure}

The streamlines are shown in Fig. (\ref{Streamline1}) for this case with
$\beta=0.5$ and in Fig. (\ref{Streamline2}) for the case $\beta=2$. For the
case $\beta=0.5$ the streamlines are topologically similar as $\alpha$ is
varied. For the case $\beta=2$ the streamlines undergo topological change as
$\alpha$ changes through $\alpha=3$ to $\alpha=5$ with the introduction of a
saddle stagnation point.%

%

\begin{figure}
[h!]
\begin{center}
\includegraphics[
natheight=5.521000in,
natwidth=4.791900in,
height=4.2047in,
width=3.653in
]%
{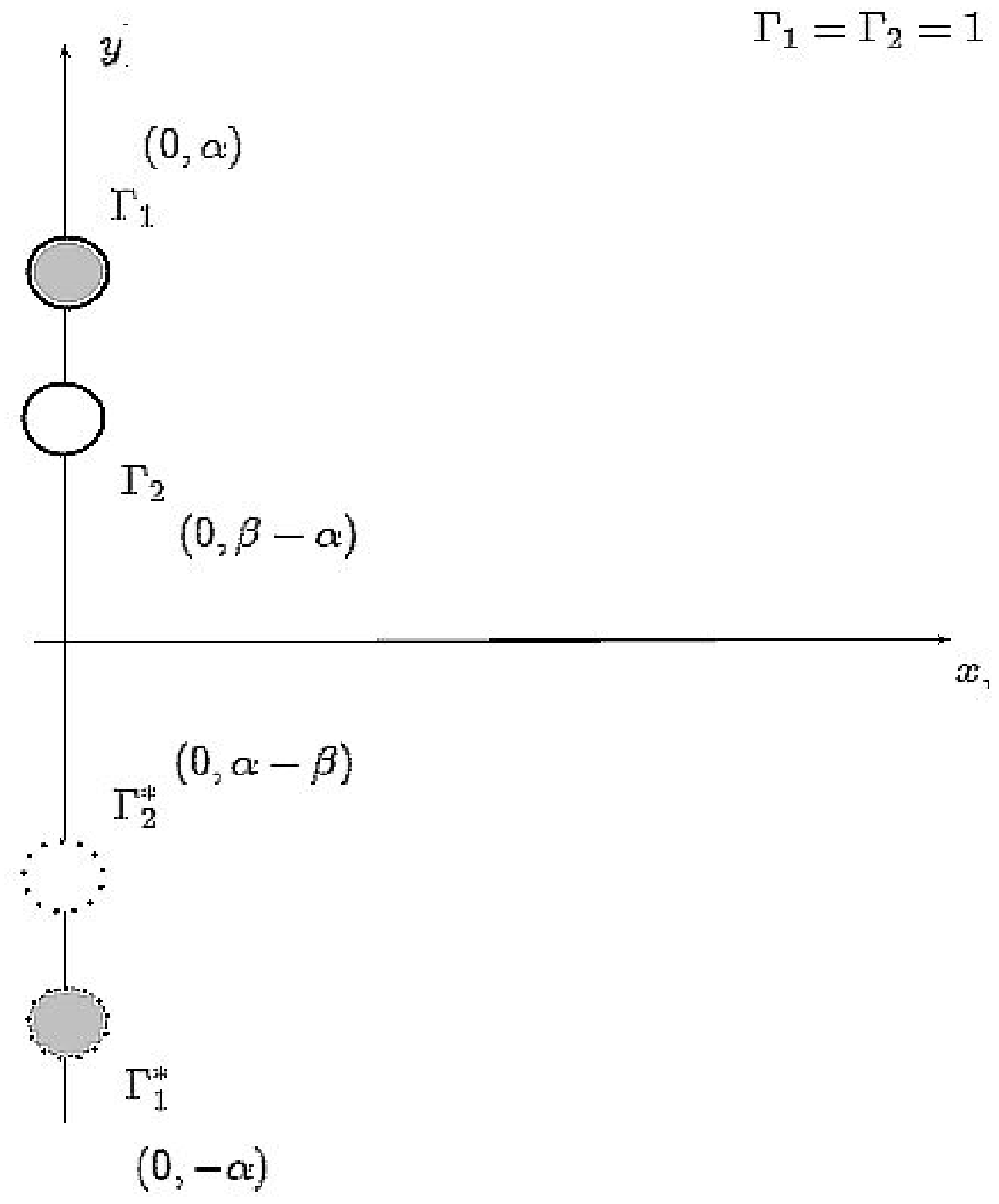}%
\caption{The configuration for streamline topologies with parameters $\alpha$,
$\beta$. Here $\Gamma_{1}$ is at $\left(  0,\alpha\right)  $, and $\Gamma_{2}$
is at $\left(  0,\beta-\alpha\right)  $ and $y_{1}+y_{2}=\beta$ is invariant.
The streamline is $\psi_{1}=\ln\sqrt{x^{2}+(y-\alpha)^{2}}-K_{0}\left(
\sqrt{x^{2}+(y-\alpha)^{2}}\right)  +K_{0}\left(  \sqrt{x^{2}+(y+\alpha)^{2}%
}\right)  -\ln\sqrt{x^{2}+(y+\alpha)^{2}}+\ln\sqrt{x^{2}+(y+\alpha-\beta)^{2}%
}+K_{0}\left(  \sqrt{x^{2}+(y+\alpha-\beta)^{2}}\right)  -K_{0}\left(
\sqrt{x^{2}+(y+\beta-\alpha)^{2}}\right)  -\ln\left(  \sqrt{x^{2}%
+(y+\beta-\alpha)^{2}}\right)  $}%
\label{Stream_configuration2}%
\end{center}
\end{figure}

The second basic configuration we consider is with 2 vortices placed in a
vertical configuration, with $\Gamma_{1}$ in layer $1$ and $\Gamma_{2}$ in
layer $2$ as shown in Fig. (\ref{Stream_configuration2}). In this case
$x_{1}=x_{2}=0$.%

\begin{figure}
[h!]
\begin{center}
\includegraphics[
natheight=5.719000in,
natwidth=5.677500in,
height=5.7467in,
width=5.7052in
]%
{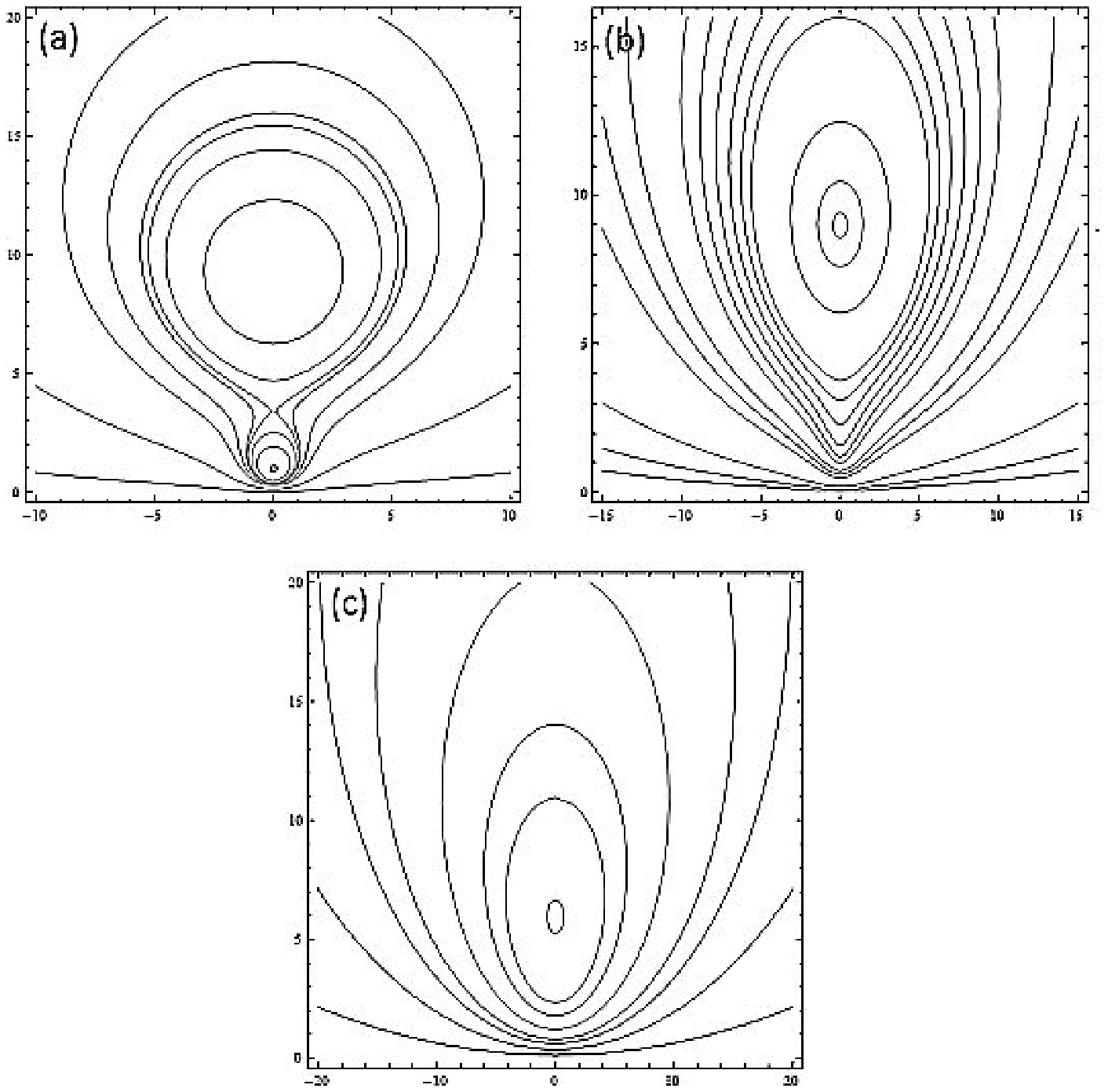}%
\caption{Streamlines for vertical configurations. $\Gamma_{1}=\Gamma_{2}$ in
layers $1$ and $2$ respectively with $\Gamma_{1}$ at $\left(  0,\alpha\right)
$ and $\Gamma_{2}$ at $\left(  0,\beta-\alpha\right)  $. In this case
$y_{1}+y_{2}=\beta$. (a) $\alpha<<\beta$; observed for $\left(  \alpha
,\beta\right)  =\left(  0.1,2\right)  ;\left(  \alpha,\beta\right)  =\left(
1,10\right)  ;\left(  \alpha,\beta\right)  =\left(  25,30\right)  $ (b)
$\alpha\approx\beta,$ observed for $\left(  \alpha,\beta\right)  =\left(
9,10\right)  ;\left(  \alpha,\beta\right)  =\left(  29,30\right)  $ (c)
$\alpha$ not much smaller than $\beta$ but not too close in magnitude to
$\beta$ either, observed for $\left(  \alpha,\beta\right)  =\left(
0.5,1\right)  ;\left(  \alpha,\beta\right)  =\left(  1,2\right)  ;\left(
\alpha,\beta\right)  =\left(  6,10\right)  ;$}%
\label{Streamline4}%
\end{center}
\end{figure}

In this case 3 classes of streamlines are observed. Notice that $\alpha<\beta
$, and $y_{1}+y_{2}=\beta$ is invariant. The first streamline case corresponds
to when $\alpha<<\beta$ (when $\Gamma_{1}$ and $\Gamma_{2\text{ }}$are
relatively far apart) and is shown in Fig. (\ref{Streamline4}(a)). There is a
saddle stagnation point observed for this case. The second case corresponds to
when $\alpha\approx\beta$ (so that is close to the boundary) and is shown in
Fig. (\ref{Streamline4}(b)). The third case corresponds to the intermediate
case when $\alpha$ is not much smaller than $\beta$ but not too similar in
magnitute to \ $\beta$ either and is depicted in Fig. (\ref{Streamline4}(c)).
For the cases $\alpha\approx\beta$ and the intermediate case when $\alpha$ is
not much smaller than $\beta$ but not too similar in magnitute to \ $\beta$
there are no stagnation points. \ While in both of these cases the streamlines
are in many ways topologically similar, the two layer interaction does
introduce distortion effects upon a closer examination of (\ref{Streamline4}%
(b)) and (\ref{Streamline4}(c)).%

\subsection{The case $\Gamma_{1}=-$ $\Gamma_{2}=1$}

\begin{figure}
[h]
\begin{center}
\includegraphics[
natheight=4.083600in,
natwidth=7.895700in,
height=2.4785in,
width=4.766in
]%
{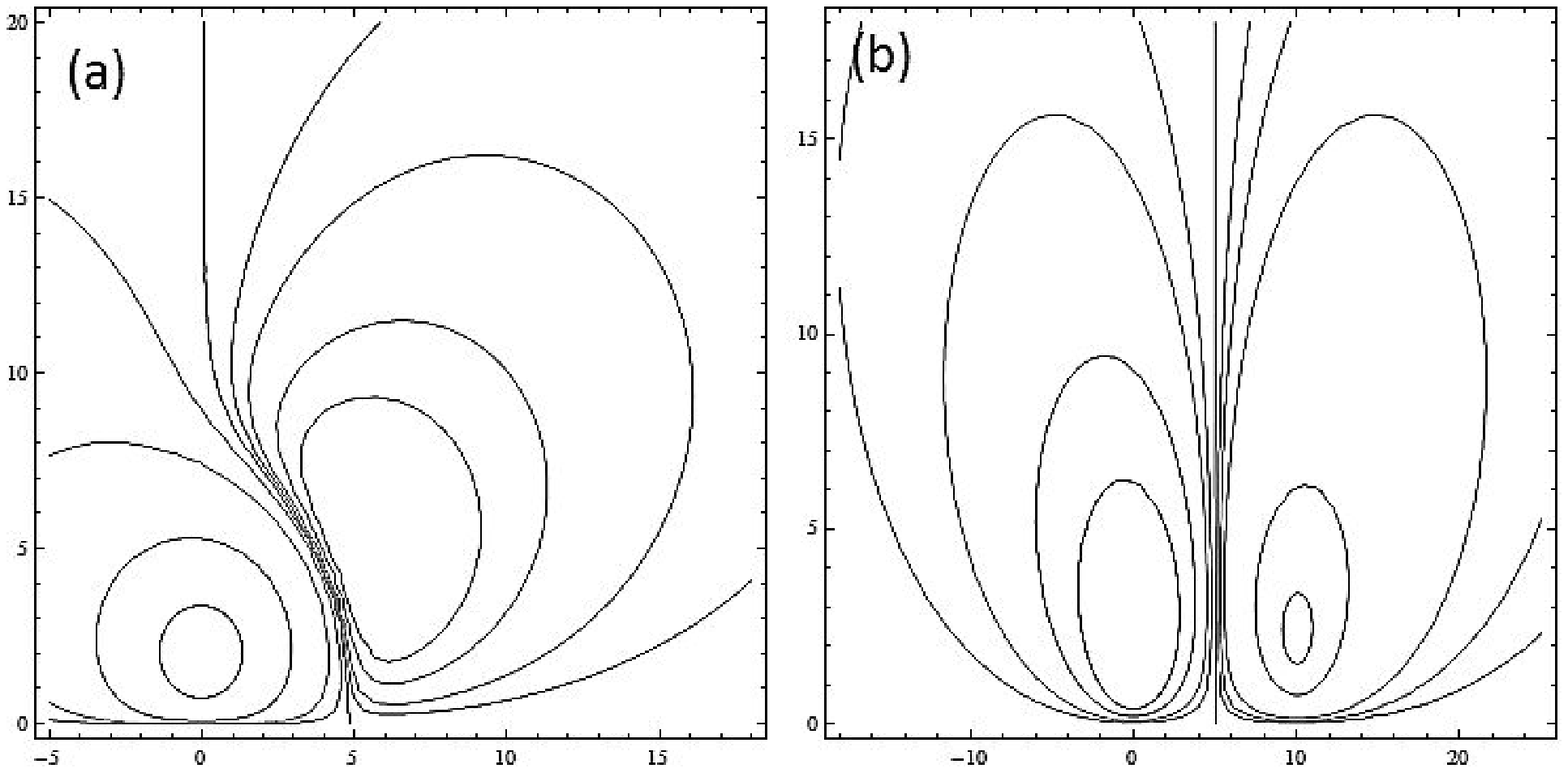}%
\caption{Streamlines for $\Gamma_{1\text{ }}$in layer 1 at $\left(
0,\beta\right)  $, $\Gamma_{2\text{ }}$in layer 2 at $\left(  \alpha
,\beta\right)  $.In this case $y_{1}-y_{2}=0$ is invariant. (a) $\beta
=2;\alpha=0.1$ and (b)$\beta=2;\alpha=10$. The streamfunction is $\psi_{1}%
=\ln\sqrt{x^{2}+(y-\beta)^{2}}-K_{0}\left(  \sqrt{x^{2}+(y-\beta)^{2}}\right)
+K_{0}\left(  \sqrt{x^{2}+(y+\beta)^{2}}\right)  -\ln\sqrt{x^{2}+(y+\beta
)^{2}}-\ln\sqrt{\left(  x-\alpha\right)  ^{2}+(y-\beta)^{2}}-K_{0}\left(
\sqrt{\left(  x-\alpha\right)  ^{2}+(y-\beta)^{2}}\right)  +K_{0}\left(
\sqrt{\left(  x-\alpha\right)  ^{2}+(y+\beta)^{2}}\right)  +\ln\left(
\sqrt{\left(  x-\alpha\right)  ^{2}+(y+\beta)^{2}}\right)  $}%
\label{Streamline5}%
\end{center}
\end{figure}

We consider the same two basic configurations. One in which $\Gamma_{1}=$
$-\Gamma_{2}=1$, with $\Gamma_{1}\ $\ in layer $1$, and $\Gamma_{2}\ $\ in
layer $2$, with the same $y$-coordinate as shown in Fig
(\ref{Stream_configuration1}). We consider various values of the parameters
$\alpha$ and $\beta$ with reference to Rossby radius of deformation $k=1$.
Note in this case $\Gamma_{1}y_{1}+\Gamma_{2}y_{2}=y_{1}-y_{2}=0$ is
invariant. \ The streamlines are shown for $\beta=2$ in Fig (\ref{Streamline5}%
). \ For all values of $\beta$ (including $\beta=0.5;2;10$) a stagnation
streamline emanating from the bounary is observed. There is a stagnation point
at the boundary from which the stagnation streamline originates. As seen for
the case $\beta=2$ in Fig (\ref{Streamline5}) the basic topology of the
streamline does not change as $\alpha$ is varied from $\alpha=0.1$ through
$\alpha=10$.%

\begin{figure}
[thb!]
\begin{center}
\includegraphics
{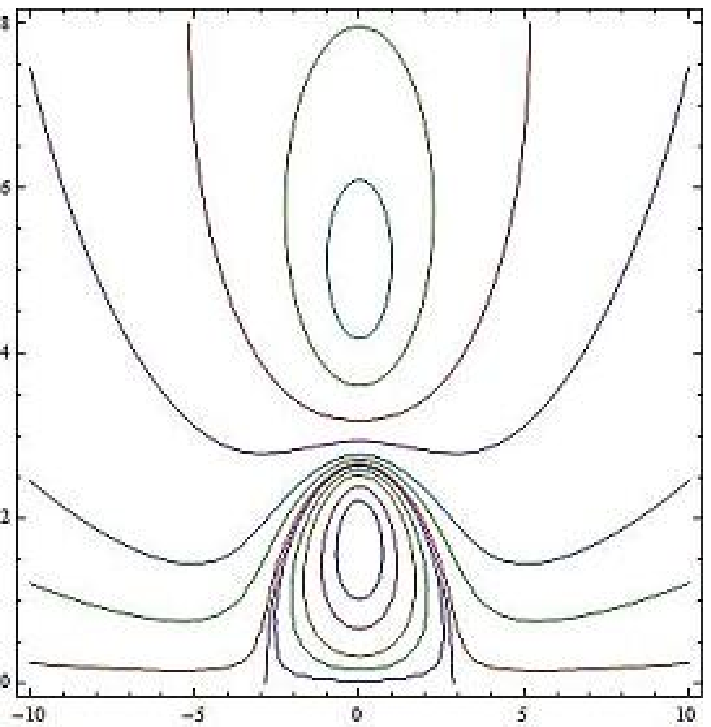}%
\caption{Streamlines for $\Gamma_{1\text{ }}$in layer 1 at $\left(
0,\alpha\right)  $, $\Gamma_{2\text{ }}$in layer 2 at $\left(  0,\alpha
-\beta\right)  $.In this case $y_{1}-y_{2}=\beta$ is invariant. The typical
case $\beta=3;\alpha=5$ is shown. The streamfunction is $\psi_{1}=\ln
\sqrt{x^{2}+(y-\alpha)^{2}}-K_{0}\left(  \sqrt{x^{2}+(y-\alpha)^{2}}\right)
+K_{0}\left(  \sqrt{x^{2}+(y+\alpha)^{2}}\right)  -\ln\sqrt{x^{2}%
+(y+\alpha)^{2}}-\ln\sqrt{x^{2}+(y-\left(  \alpha-\beta\right)  )^{2}}%
-K_{0}\left(  \sqrt{x^{2}+(y-\left(  \alpha-\beta\right)  )^{2}}\right)
+K_{0}\left(  \sqrt{x^{2}+(y+\left(  \alpha-\beta\right)  )^{2}}\right)
+\ln\left(  \sqrt{x^{2}+(y+\left(  \alpha-\beta\right)  )^{2}}\right)$}
\label{Streamline6}%
\end{center}
\end{figure}
The second configuration is as before with $\Gamma_{1}=$ $-\Gamma_{2}=1$, with
$\Gamma_{1}\ $\ in layer $1$ and $\Gamma_{2}\ $\ in layer $2$ . The 2 vortices
are placed in a vertical configuration, with $\Gamma_{1}$ in layer $1$ at
$\left(  0,\alpha\right)  $ and $\Gamma_{2}$ in layer $2$ at $\left(
0,\alpha-\beta\right)  $ shown in Fig. (\ref{Stream_configuration2}). In this
case $x_{1}=x_{2}=0$. We consider various values of the parameters $\alpha$
and $\beta$ with reference to Rossby radius of deformation $k=1$. Note in this
case $\Gamma_{1}y_{1}+\Gamma_{2}y_{2}=y_{1}-y_{2}=\beta$ is invariant. \ In
this case there is a stagnation streamline joining two stagnation points on
the boundary as seen in Fig. (\ref{Streamline6}).%
\section{Integrable 3 vortex configurations--Relative equilibria.}

\noindent We conclude with a Hamiltonian study of integrable two-layer 3
vortex motion on the half plane by studying integrable symmetrical
configurations and provide a rich class of new relative equilibria. We
consider two basic symmetrical configurations of 3 vortices as depicted in
Fig.(\ref{Int_3_vortex_config}). We seek relative equilibrium solutions.in
which the initial configuration is rigidly maintained. \ The method we adopt
is similar to that as in the work of Jamaloodeen and Newton \cite{JamNewt2} in
which we seek relative equilbria base on a a symmetrical configuration and
vary parameters that ensure the relative equilbrium or rigid shape is invariant.%

\subsection{The case $\Gamma_{1}=$ $\Gamma_{2}=1,\Gamma_{3}=-\alpha$ with
$\Gamma_{1}$ at $\left(  x,y\right)  $, $\Gamma_{2}$ at $\left(  -x,y\right)
$ (in layer 1) and $\Gamma_{3}$ at $\left(  0,y_{3}\right)  $ (in layer 2)}
\begin{figure}
[h]
\begin{center}
\includegraphics[
natheight=6.146200in,
natwidth=9.729100in,
height=3.1401in,
width=4.9588in
]%
{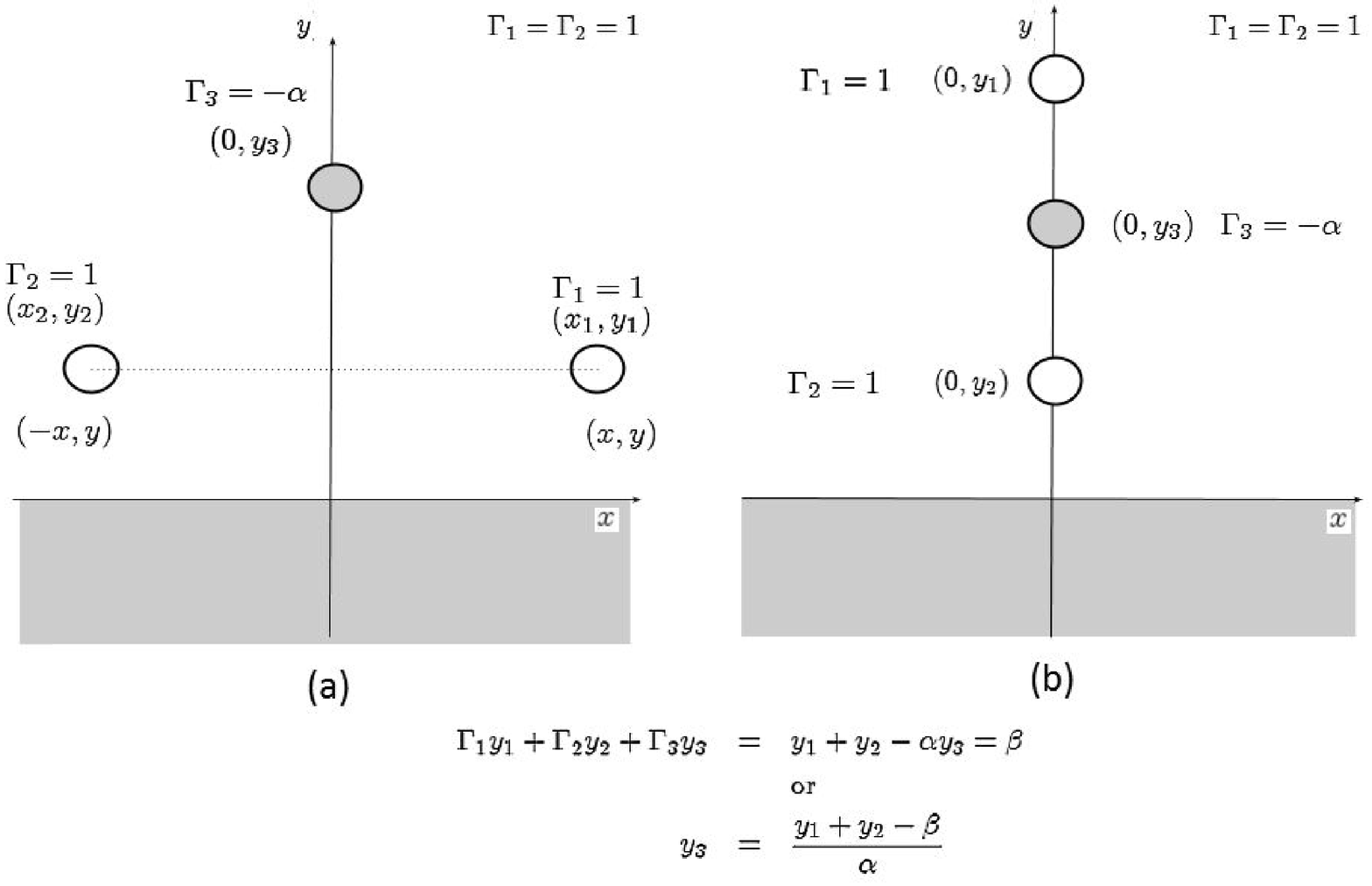}%
\caption{Integrable 3 vortex relative equilibrium solutions depending on
parameters, $\alpha$, and $\beta$. The invariant $\Gamma_{1}y_{1}+\Gamma
_{2}y_{2}+\Gamma_{3}y_{3}=\beta$ becomes $y_{1}+y_{2}-\alpha y_{3}=\beta$ with
$\Gamma_{1}=\Gamma_{2}=1$, $\Gamma_{3}=-\alpha$, so that $y_{3}=\frac
{y_{1}+y_{2}-\beta}{\alpha}$. $\Gamma_{1}$ and $\Gamma_{2}$ are in layer 1,
$\Gamma_{3}$ is in layer 2. In (a) $\Gamma_{1}$ and $\Gamma_{2}$ are placed at
$\left(  x,y\right)  $ and $\left(  -x,y\right)  $ respectively. In (b)
$\Gamma_{1}$ and $\Gamma_{2}$ are placed at $\left(  0,y_{1}\right)  $ and
$\left(  0,y_{2}\right)  $ respectively.}%
\label{Int_3_vortex_config}%
\end{center}
\end{figure}

\noindent In the first symmetrical configuration we have $\Gamma_{1}=$ $\Gamma_{2}=1$,
with $\Gamma_{1}\ $\ in layer $1$, and $\Gamma_{2}\ $\ in layer $2$, with the
same $y$-coordinate as shown in Fig (\ref{Int_3_vortex_config})(a). In this
case it can be shown that $\dot{x}_{1}=\dot{x}_{2}$ and $\dot{y}_{1}=-\dot
{y}_{2}$. The invariant $\Gamma_{1}y_{1}+\Gamma_{2}y_{2}+\Gamma_{3}y_{3}%
=\beta$, becomes $y_{1}+y_{2}-\alpha y_{3}=\beta$, or $y_{3}=\frac{y_{1}%
+y_{2}-\beta}{\alpha}$. Relative equilibrium solutions will then be admitted,
for this configuration, by requiring that $\dot{y}_{1}=0$, and $\dot{x}%
_{1}=\dot{x}_{3}$. These equations (see Appendix B for their explicit forms)
are solved numerically with $\alpha$ and $\beta$ considered as parameters.
These are summarized in in Tables (\ref{Horizontal_3_equilibria_beta_negative}%
)-(\ref{Horizontal_3_equilibria_beta_positive}), and exploiting the symmetry
and the invariant $y_{1}+y_{2}-\alpha y_{3}=\beta$, or $y_{3}=\frac
{y_{1}+y_{2}-\beta}{\alpha},$ only the coordinates $\left(  x_{1},y_{1}\right) =\left(x,y\right)$ are given.%
\begin{equation}%
\begin{tabular}
[c]{l|l}\hline\hline
$\left(  \alpha,\beta\right)  $ (negative $\beta)$ & Relative equilibrium
$\left(  x,y\right)  $\\\hline\hline
$\left(  2.005,-0.25\right)  $ & $\left(  0.5639,4.524\right)  $\\\hline
$\left(  2.01,-0.25\right)  $ & $\left(  0.6827,4.328\right)  $\\\hline
$\left(  2.02,-0.25\right)  $ & $\left(  0.7252,4.114\right)  ,\left(
0.8457,12.43\right)  $\\\hline
$\left(  2.05,-0.5\right)  $ & $\left(  0.5847,3.291\right)  ,\left(
0.8497,9.977\right)  ,\left(  0.01254,9.941\right)  $\\\hline
$\left(  2.1,-0.5\right)  $ & $\left(  0.8317,4.851\right)  ,\left(
0.7775,3.466\right)  $\\\hline
$\left(  2.11,-0.5\right)  $ & $\left(  0.8558,3.857\right)  $\\\hline
$\left(  2.01,-0.75\right)  $ & $\left(  0.4542,3.737\right)  $\\\hline
$\left(  2.09,-0.75\right)  $ & $\left(  0.8156,8.314\right)  ,\left(
0.4843,2.773\right)  $\\\hline
$\left(  2.15,-0.75\right)  $ & $\left(  0.7855,4.941\right)  ,\left(
0.665,2.863\right)  $\\\hline
$\left(  2.19,-0.75\right)  $ & $\left(  0.8015,3.917\right)  $\\\hline
$\left(  2.1,-1\right)  $ & $\left(  0.5144,2.569\right)  ,\left(
0.8156,10.01\right)  $\\\hline
$\left(  2.25,-1\right)  $ & $\left(  0.7955,3.857\right)  ,\left(
0.663,2.622\right)  $\\\hline
$\left(  2.27,-1\right)  $ & $\left(  0.7775,3.315\right)  ,\left(
0.7172,2.773\right)  $\\\hline
$\left(  2.275,-1\right)  $ & $\left(  0.7654,3.104\right)  $\\\hline
$\left(  2.7,-5\right)  $ & $\left(  0.7252,7.11\right)  ,\left(
0.9661,1.177\right)  $\\\hline
$\left(  3,-5\right)  $ & $\left(  0.916,0.9963\right)  ,\left(
0.657,5.002\right)  ,\left(  0.04868,4.941\right)  $\\\hline
$\left(  4,-5\right)  $ & $\left(  0.5426,2.442\right)  ,\left(
0.7232,0.9059\right)  ,\left(  0.06072,2.442\right)  $\\\hline
$\left(  4.5,-5\right)  $ & $\left(  0.4522,1.779\right)  ,\left(
0.6028,0.9963\right)  $\\\hline
$\left(  4.825,-5\right)  $ & $\left(  0.4348,1.313\right)  $\\\hline
\end{tabular}
\label{Horizontal_3_equilibria_beta_negative}%
\end{equation}
\begin{equation}%
\begin{tabular}
[c]{l|l}\hline\hline
$\left(  \alpha,\beta\right)  $ (positive $\beta)$ & Relative equilibrium
$\left(  x,y\right)  $\\\hline\hline
$\left(  1.9855,0.1\right)  $ & $\left(  0.8845,6.198\right)  $\\\hline
$\left(  1.99,0.1\right)  $ & $\left(  0.8678,9.977\right)  ,\left(
0.8497,5.46\right)  $\\\hline
$\left(  1.9755,0.15\right)  $ & $\left(  0.8964,5.509\right)  $\\\hline
$\left(  1.98,0.15\right)  $ & $\left(  0.8608,7.385\right)  ,\left(
0.8548,4.987\right)  $\\\hline
$\left(  1.9533,0.25\right)  $ & $\left(  0.8798,4.791\right)  $\\\hline
$\left(  1.954,0.25\right)  $ & $\left(  0.8678,5.152\right)  ,\left(
0.8738,4.58\right)  $\\\hline
$\left(  (1.96,0.25\right)  $ & $\left(  0.8738,6.206\right)  ,\left(
0.8558,4.399\right)  $\\\hline
$\left(  1.888,0.5\right)  $ & $\left(  0.91,3.917\right)  $\\\hline
$\left(  1.9,0.5\right)  $ & $\left(  0.8919,4.917\right)  ,\left(
0.8317,3.544\right)  $\\\hline
$\left(  1.95,0.5\right)  $ & $\left(  0.5907,3.448\right)  ,\left(
0.8738,10.03\right)  ,(0.01856,10.03)$\\\hline
$\left(  1.807,0.75\right)  $ & $\left(  0.9201,3.491\right)  $\\\hline
$\left(  1.82,0.75\right)  $ & $\left(  0.9201,4.052\right)  ,\left(
0.8608,3.158\right)  $\\\hline
$\left(  1.715,1\right)  $ & $\left(  0.9558,3.075\right)  $\\\hline
$\left(  1.75,1\right)  $ & $\left(  0.932,3.927\right)  ,\left(
0.8192,2.763\right)  $\\\hline
$\left(  1.82,1\right)  $ & $\left(  0.9201,5.569\right)  ,\left(
0.653,2.743\right)  $\\\hline
$\left(  1.9,1\right)  $ & $\left(  0.5546,2.93\right)  ,\left(
0.8919,9.941\right)  ,\left(  0.01254,9.941\right)  $\\\hline
$\left(  1.483,1.5\right)  $ & $\left(  1.009,2.541\right)  $\\\hline
$\left(  1.5,1.5\right)  $ & $\left(  1.027,2.802\right)  ,\left(
0.9617,2.375\right)  $\\\hline
$\left(  1.7,1.5\right)  $ & $\left(  0.9616,5.022\right)  ,\left(
0.7479,2.44\right)  $\\\hline
$\left(  1.85,1.5\right)  $ & $\left(  0.9039,10.01\right)  ,\left(
0.01254,10.01\right)  ,\left(  0.6871,2.858\right)  $\\\hline
$\left(  1.34,1.75\right)  $ & $\left(  1.07,2.334\right)  $\\\hline
$\left(  1.4,1.75\right)  $ & $\left(  0.9912,2.179\right)  ,\left(
1.082,2.79\right)  $\\\hline
$\left(  1.18,2\right)  $ & $\left(  1.181,2.161\right)  $\\\hline
$\left(  1.25,2\right)  $ & $\left(  1.098,2.048\right)  ,\left(
1.164,2.553\right)  $\\\hline
$\left(  1.5,2\right)  $ & $\left(  1.045,4.007\right)  ,\left(
0.9379,2.262\right)  $\\\hline
$\left(  0.8,2.5\right)  $ & $\left(  1.532,1.888\right)  $\\\hline
$\left(  0.85,2.5\right)  $ & $\left(  1.496,1.852\right)  ,\left(
1.472,2.066\right)  $\\\hline
$\left(  1.5,2.5\right)  $ & $\left(  1.051,4.987\right)  ,\left(
1.042,2.44\right)  $\\\hline
$\left(  1.75,2.5\right)  $ & $\left(  0.9494,3.03\right)  ,\left(
0.9403,10.03\right)  ,\left(  0.02784,10.03\right)  $\\\hline
\end{tabular}
\ \label{Horizontal_3_equilibria_beta_positive}%
\end{equation}

\begin{figure}
[h!]
\begin{center}
\includegraphics
{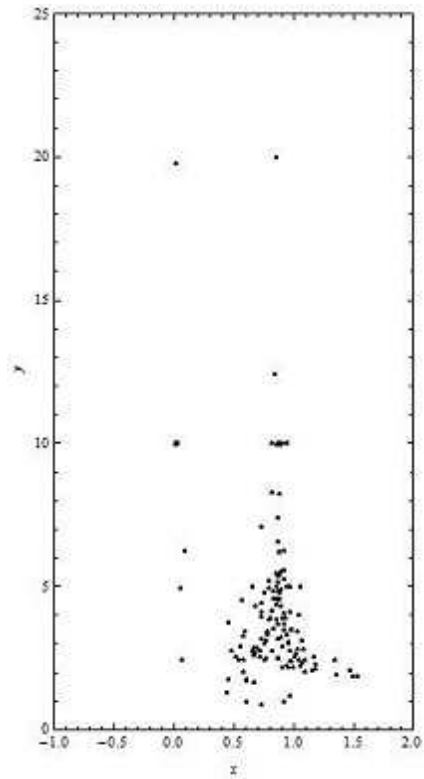}%
\caption{A scatterplot of all numerically found relative equilibria from
Tables (\ref{Horizontal_3_equilibria_beta_negative}%
)-(\ref{Horizontal_3_equilibria_beta_positive}) with associated $\left(
\alpha,\beta\right)  $ for each configuration $\left(  x_{1},y_{1}\right)
=\left(  x,y\right)  $ suppressed. }%
\label{s_plot1}%
\end{center}
\end{figure}

Notice that the numerical evidence suggest that there are no relative
equilibria corresponding to this configuration with $\alpha$ negative or
$\Gamma_{3}>0$. There also appears to be a complicated bifurcation of these
relative equilibria. For example fixing $\beta$ and varying $\alpha$ gives
varying numbers of relative equilibria. Consider, for example the case
$\beta=-0.75$ for which at $\alpha=2.01$, there is one relative equilibrium
solution, then increasing $\alpha$ through $\alpha=2.09$ there are two and
likewise two, again, after increasing $\alpha$ through $\alpha=2.15$ and
finally increasing $\alpha$ to $\alpha=2.19$ there is again only one relative
equilibrium solution. Similarly when $\beta=0.5$ there is one equilibrium
configuration when $\alpha=1.888$, two equilibrium configurations when
$\alpha=1.9$ and three equilibrium configurations when $\alpha$ increases
through $\alpha=1.95$. A scatterplot of all numerically found relative
equilibria from Tables (\ref{Horizontal_3_equilibria_beta_negative}%
)-(\ref{Horizontal_3_equilibria_beta_positive}) is shown in Fig(\ref{s_plot1}%
), with associated $\left(  \alpha,\beta\right)  $ for each configuration
$\left(  x_{1},y_{1}\right)  =\left(  x,y\right)  $ suppressed.

It is also possible that some of these relative equilibria are fixed
equilibria, meaning that in fact not only are $\dot{y}_{1}=0$, and $\dot
{x}_{1}=\dot{x}_{3}$ but $\dot{y}_{1}=0$, and $\dot{x}_{1}=\dot{x}_{3}=0$.
However preliminary numerical studies suggest that there are no fixed
equilibrium configurations of the kind depicted in
Fig(\ref{Int_3_vortex_config})(a). \ The bifurcation sequence suggested above
and the fixed equilibria are topics for further study.

\subsection{The case $\Gamma_{1}=$ $\Gamma_{2}=1,\Gamma_{3}=-\alpha$ with
$\Gamma_{1}$ at $\left(  0,y_{1}\right)  $, $\Gamma_{2}$ at $\left(
0,y_{2}\right)  $ (in layer 1) and $\Gamma_{3}$ at $\left(  0,y_{3}\right)  $
(in layer 2)}

In the second symmetrical configuration we have $\Gamma_{1}=$ $\Gamma_{2}=1$,
with $\Gamma_{1}\ $\ in layer $1$, and $\Gamma_{2}\ $\ in layer $2$, with
$\Gamma_{1}$ at $\left(  0,y_{1}\right)  $, $\Gamma_{2}$ at $\left(
0,y_{2}\right)  $ (in layer 1) and $\Gamma_{3}$ at $\left(  0,y_{3}\right)  $
(in layer 2)with the same $y$-coordinate as shown in Fig
(\ref{Int_3_vortex_config})(b). In this case it can be shown that $\dot{x}%
_{1}=\dot{x}_{2}$ and $\dot{y}_{1}=-\dot{y}_{2}$. The invariant $\Gamma
_{1}y_{1}+\Gamma_{2}y_{2}+\Gamma_{3}y_{3}=\beta$, becomes $y_{1}+y_{2}-\alpha
y_{3}=\beta$, or $y_{3}=\frac{y_{1}+y_{2}-\beta}{\alpha}$. In this case it can
be shown that $\dot{y}_{1}=\dot{y}_{2}=\dot{y}_{3}=0$. Relative equilibrium
solutions will then be admitted, for this configuration by requiring that
$\dot{x}_{1}=\dot{x}_{2}=\dot{x}_{3}$. Again the invariant $\Gamma_{1}%
y_{1}+\Gamma_{2}y_{2}+\Gamma_{3}y_{3}=\beta$, becomes $y_{1}+y_{2}-\alpha
y_{3}=\beta$, or $y_{3}=\frac{y_{1}+y_{2}-\beta}{\alpha}$. These equations
(see Appendix B for their explicit forms) are solved numerically with $\alpha$
and $\beta$ considered as parameters. These are summarized in in Tables
(\ref{vertical_3_equilibria_beta_1})-(\ref{vertical_3_equilibria_beta_2}), and
exploiting the symmetry and the invariant $y_{1}+y_{2}-\alpha y_{3}=\beta$, or
$y_{3}=\frac{y_{1}+y_{2}-\beta}{\alpha},$ again it suffices to provide only
the coordinates for $\left(  y_{1},y_{2}\right)  =(Y,y)$.%

\begin{equation}%
\begin{tabular}
[c]{l|l}\hline\hline
$\left(  \alpha,\beta\right)  $ ($\alpha\leq2$) & Relative equilibrium
$\left(  Y,y\right)  $\\\hline\hline
$\alpha\leq1$ & No numerical relative equilibrium configurations\\\hline
$\left(  1.5,0.25\right)  $ & $\left(  4.132,1.549\right)  $\\\hline
$\left(  1.5,0.5\right)  $ & $\left(  3.859,1.306\right)  $\\\hline
$\left(  1.5,1\right)  $ & $\left(  3.687,1.211\right)  $\\\hline
$\left(  1.5,2\right)  $ & $\left(  3.758,1.353\right)  $\\\hline
$\left(  1.5,3\right)  $ & $\left(  4.369,1.749\right)  $\\\hline
$\left(  1.5,3.4\right)  $ & $\left(  4.72,1.923\right)  ,\left(
2.333,4.987\right)  $\\\hline
$\left(  1.5,3.5\right)  $ & $\left(  4.827,1.995\right)  ,\left(
2.832,5.379\right)  ,\left(  2.03,4.844\right)  $\\\hline
$\left(  1.5,5\right)  $ & $\left(  1.809,5.995\right)  ,\left(
5.393,8.073\right)  ,\left(  6.598,2.893\right)  $\\\hline
$\left(  1.5,10\right)  $ & $\left(  13.32,6.533\right)  ,\left(
13.92,17.06\right)  ,\left(  2.414,11.23\right)  $\\\hline
$\left(  1.75,-5\right)  $ & $\left(  1.037,8.751\right)  $\\\hline
$\left(  1.75,-2\right)  $ & $\left(  0.9649,4.791\right)  $\\\hline
$\left(  1.75,-1\right)  $ & $\left(  1.336,3.903\right)  $\\\hline
$\left(  1.75,-0.5\right)  $ & $\left(  1.899,4.102\right)  $\\\hline
$\left(  1.75,0.5\right)  $ & $\left(  3.68,1.736\right)  $\\\hline
$\left(  1.75,1\right)  $ & $\left(  3.703,1.745\right)  $\\\hline
$\left(  1.75,1.5\right)  $ & $\left(  3.963,1.959\right)  ,\left(
2.716,4.653\right)  $\\\hline
$\left(  1.75,1.75\right)  $ & $\left(  4.187,2.152\right)  ,\left(
1.846,3.936\right)  ,\left(  3.772,5.704\right)  $\\\hline
$\left(  1.75,2\right)  $ & $\left(  4.52,2.412\right)  ,\left(
1.659,3.857\right)  ,\left(  4.7,6.628\}\right)  $\\\hline
$\left(  1.75,3\right)  $ & $\left(  6.222,3.821\right)  ,\left(
1.436,4.303\right)  ,\left(  7.553,9.555\right)  $\\\hline
$\left(  2,0.25\right)  $ & $\left(  4.38,2.648\right)  ,\left(
2.365,4.107\right)  $\\\hline
$\left(  2,0.5\right)  $ & $\left(  4.589,2.837\right)  ,\left(
2.015,3.771\right)  $\\\hline
$\left(  2,0.75\right)  $ & $\left(  6.135,4.349\right)  ,\left(
1.745,3.562\right)  $\\\hline
$\left(  2,1\right)  $ & $\left(  1.603,3.491\right)  $\\\hline
$\left(  2,2\right)  $ & $\left(  1.3,3.651\right)  $\\\hline
$\left(  2,3\right)  $ & $\left(  1.235,4.349\right)  $\\\hline
$\left(  2,5\right)  $ & $\left(  1.302,6.296\right)  $\\\hline
$\left(  2,10\right)  $ & $\left(  1.708,11.66\right)  $\\\hline
$\left(  2,15\right)  $ & $\left(  2.201,17.14\right)  $\\\hline
\end{tabular}
\ \ \label{vertical_3_equilibria_beta_1}%
\end{equation}

\begin{equation}%
\begin{tabular}
[c]{l|l}\hline\hline
$\left(  \alpha,\beta\right)  $ ($\alpha>2$) & Relative equilibrium $\left(
Y,y\right)  $\\\hline
$\left(  2.35,-6\right)  $ & $\left(  14.69,13.01\right)  ,\left(
6.292,9.165\right)  ,\left(  0.5845,6.417\right)  $\\\hline
$\left(  2.35,-5\right)  $ & $\left(  12.34,10.67\right)  ,\left(
4.927,7.412\right)  ,\left(  0.578,5.416\right)  $\\\hline
$\left(  2.35,-3\right)  $ & $\left(  7.983,6.327\right)  ,\left(
0.665,3.405\right)  ,\left(  2.472,4.309\right)  $\\\hline
$\left(  2.35,-2\right)  $ & $\left(  5.784,4.219\right)  $\\\hline
$\left(  2.35,-1\right)  $ & $\left(  4.31,2.76\right)  $\\\hline
$\left(  2.35,-0.5\right)  $ & $\left(  4.615,3.011\right)  $\\\hline
$\left(  2.35,0.25\right)  $ & $\left(  2.012,3.651\right)  $\\\hline
$\left(  2.35,1\right)  $ & $(1.365,3.28)$\\\hline
$\left(  2.35,1.5\right)  $ & $\left(  1.223,3.408\right)  $\\\hline
$\left(  2.35,1.75\right)  $ & $\left(  1.202,3.547\right)  $\\\hline
$\left(  2.35,2\right)  $ & $\left(  1.175,3.704\right)  $\\\hline
$\left(  2.35,5\right)  $ & $\left(  1.211,6.567\right)  $\\\hline
$\left(  2.35,10\right)  $ & $\left(  1.538,11.87\right)  $\\\hline
$\left(  3,-10\right)  $ & $\left(  9.041,7.594\right)  ,\left(
32.93,21.49\}\right)  ,\left(  3.121,6.554\right)  $\\\hline
$\left(  3,-5\right)  $ & $\left(  4.898,3.509\right)  ,\left(
16.3,10.69\right)  ,\left(  1.246,3.206\right)  $\\\hline
$\left(  3,-4\right)  $ & $\left(  4.114,2.755\right)  ,\left(
12.78,8.437\right)  $\\\hline
$\left(  3,-3\right)  $ & $\left(  3.441,2.125\right)  ,\left(
8.962,6.012\right)  $\\\hline
$\left(  3,-2\right)  $ & $\left(  3.341,2.003\right)  ,\left(
5.525,3.685\right)  $\\\hline
$\left(  3,0.25\right)  $ & $\left(  1.852,3.954\right)  $\\\hline
$\left(  3,1\right)  $ & $\left(  1.247,3.488\right)  $\\\hline
$\left(  3,2\right)  $ & $\left(  1.086,4.096\right)  $\\\hline
$\left(  3,5\right)  $ & $\left(  1.131,7.039\right)  $\\\hline
$\left(  3,10\right)  $ & $\left(  1.387,12.53\right)  $\\\hline
$\left(  5,-10\right)  $ & $\left(  3.396,2.36\right)  ,\left(
6.055,4.096\right)  $\\\hline
$\left(  5,0.25\right)  $ & $\left(  4.301,17.17\right)  $\\\hline
$\left(  5,0.5\right)  $ & $\left(  2.729,10.86\right)  $\\\hline
$\left(  5,1\right)  $ & $\left(  1.796,7.29\right)  $\\\hline
$\left(  5,5\right)  $ & $\left(  1.176,8.581\right)  $\\\hline
$\left(  5,5\right)  $ & $\left(  1.248,13.86\right)  $\\\hline
\end{tabular}
\ \ \ \label{vertical_3_equilibria_beta_2}%
\end{equation}

\bigskip Notice again that the numerical evidence suggests that there are no
relative equilibria corresponding to this configuration with $\alpha$ negative
or $\Gamma_{3}>0$. There again appears to be a complicated bifurcation of
these relative equilibria. For example fixing $\alpha$ and varying $\beta$
gives varying numbers of relative equilibria. A scatterplot of all numerically
found relative equilibria from Tables (\ref{vertical_3_equilibria_beta_1}%
)-(\ref{vertical_3_equilibria_beta_2}) is shown in Fig(\ref{s_plot1}), with
associated $\left(  \alpha,\beta\right)  $ for each configuration $\left(
0,y_{1}\right)  =\left(  0,Y\right)  \ $and $\left(  0,y_{2}\right)  =\left(
0,y\right)  $ suppressed.%
\begin{figure}
[h!]
\begin{center}
\includegraphics
{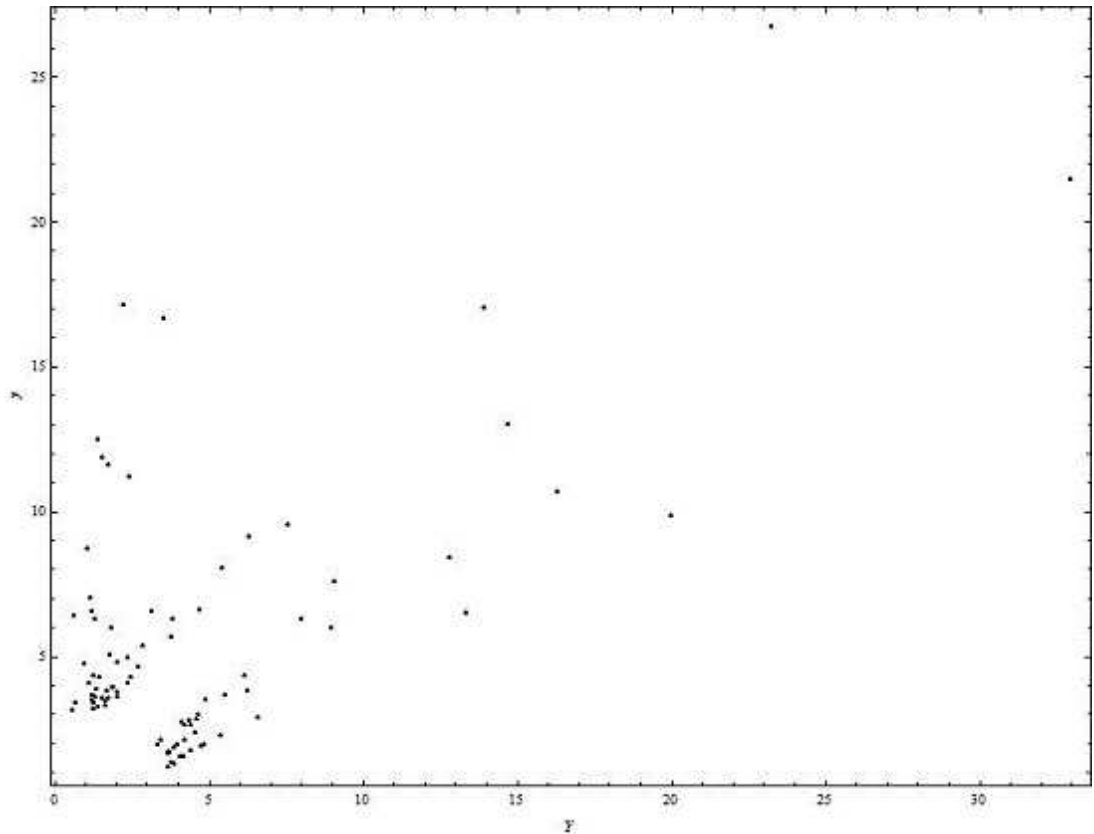}%
\caption{A scatterplot of all numerically found relative equilibria from
Tables (\ref{vertical_3_equilibria_beta_1}%
)-(\ref{vertical_3_equilibria_beta_2}) with associated $\left(  \alpha
,\beta\right)  $ for each configuration $\left(  0,y_{1}\right)  =\left(
0,Y\right)  \ $and $\left(  0,y_{2}\right)  =\left(  0,y\right)  $ or $\left(
Y,y\right)  $ suppressed. }%
\label{s_plot2}%
\end{center}
\end{figure}

\section{\bigskip Conclusion}

We have presented results on integrable two-layer point vortex motion on the
upper half plane and shown similarities and differences with integrable
\ one-layer point vortex motion in the upper half plane and with integrable
two-layer point vortex motion on the entire plane including the study of
equilbrium solutions, the finite-time vortex collapse problem and a study of
streamline topologies. \ At present we are also pursuing work to determine 3
or 4 vortex initial vortex configurations that lead to finite time collapse of
the vortices. Should such configurations exist, they will depend on the
initial configuration of the vortices and the vortex strengths, and would
possibly be self-similar. Also of interest would be to perform a systematic
bifurcation analysis of the 3 vortex relative equilibria we found shown in
Figs.(\ref{s_plot1}-\ref{s_plot2}). The bifurcation analysis in both cases
would depend on the parameter $\left(  \alpha,\beta\right)  $ as depicted in
Fig. (\ref{Int_3_vortex_config}). Finally we are also pursuing work in the
direction of establishing conditions under which 3 vortices may induce chaotic
advection of fluid particles much along the lines for the two layer problem in
the entire plane as done in the work of Koshel et al\cite{KSV2}. A good
starting point are perturbations of the configurations used to obtain the
streamlines depicted in Figs. (\ref{Stream_configuration1} and
\ref{Stream_configuration2}).

\appendix

\section{Explicit proof that $\sum\Gamma_{i}=c$ is invariant}

We show that $\Gamma_{1}\dot{y}_{1}+\Gamma_{2}\dot{y}_{2}=0$. Use,

\bigskip%
\begin{equation}
\dot{x}_{i}=-\left.  \frac{\partial\psi_{i}}{\partial y}\right\vert _{\left(
x_{i},y_{i}\right)  },\text{ \ \ \ \ \ \ }\dot{y}_{i}=\left.  \frac
{\partial\psi_{i}}{\partial x}\right\vert _{_{\left(  x_{i},y_{i}\right)  }%
}\text{ \ }%
\end{equation}

\begin{multline}
\text{\ }\dot{y}_{1}=\left.  \frac{\partial\psi_{1}}{\partial x}\right\vert
_{_{\left(  x_{1},y_{1}\right)  }}=\Gamma_{1}\left(  -K_{1}\left(
2y_{1}\right)  \right)  \frac{2\left(  x_{1}-x_{1}\right)  }{2r_{1,1^{\ast}}%
}-\frac{\Gamma_{1}}{2r_{1,1^{\ast}}^{2}}2\left(  x_{1}-x_{1}\right)
+\frac{\Gamma_{2}}{2r_{1,2}^{2}}2\left(  x_{1}-x_{2}\right) \\
-\Gamma_{2}K_{1}\left(  r_{1,2}\right)  \frac{2\left(  x_{1}-x_{2}\right)
}{2r_{1,2}}+\Gamma_{2}K_{1}\left(  r_{1,2^{\ast}}\right)  \frac{2\left(
x_{1}-x_{2}\right)  }{2r_{1,2^{\ast}}}-\frac{\Gamma_{2}}{2r_{1,2^{\ast}}^{2}%
}2\left(  x_{1}-x_{2}\right) \\
=\Gamma_{2}\left(  x_{1}-x_{2}\right)  \left[  \frac{1}{r_{1,2}^{2}}%
-\frac{K_{1}\left(  r_{1,2}\right)  }{r_{1,2}}+\frac{K_{1}\left(
r_{1,2^{\ast}}\right)  }{r_{1,2^{\ast}}}-\frac{1}{r_{1,2^{\ast}}^{2}}\right]
\end{multline}

\begin{multline}
\text{\ }\dot{y}_{2}=\left.  \frac{\partial\psi_{1}}{\partial x}\right\vert
_{_{\left(  x_{2},y_{2}\right)  }}=\Gamma_{1}K_{1}\left(  2y_{1}\right)
\frac{2\left(  x_{2}-x_{2}\right)  }{2r_{2,2^{\ast}}}-\frac{\Gamma_{2}%
}{2r_{2,2^{\ast}}^{2}}2\left(  x_{2}-x_{2}\right)  +\frac{\Gamma_{1}}%
{2r_{1,1}^{2}}2\left(  x_{2}-x_{1}\right) \\
-\Gamma_{1}K_{1}\left(  r_{2,1}\right)  \frac{2\left(  x_{2}-x_{1}\right)
}{2r_{2,1}}+\Gamma_{1}K_{1}\left(  r_{2,1^{\ast}}\right)  \frac{2\left(
x_{2}-x_{1}\right)  }{2r_{2,1^{\ast}}}-\frac{\Gamma_{1}}{2r_{2,1^{\ast}}^{2}%
}2\left(  x_{2}-x_{1}\right) \\
=\Gamma_{1}\left(  x_{2}-x_{1}\right)  \left[  \frac{1}{r_{2,1}^{2}}%
-\frac{K_{1}\left(  r_{2,1}\right)  }{r_{2,1}}+\frac{K_{1}\left(
r_{2,1^{\ast}}\right)  }{r_{2,1^{\ast}}}-\frac{1}{r_{2,1^{\ast}}^{2}}\right]
\text{.}%
\end{multline}

Here $K_{0}^{^{\prime}}(r)=-K_{1}\left(  r\right)  $ It is clear that since
$r_{2,1}=r_{2,1}$ and $r_{1,2^{\ast}}=r_{2,1^{\ast}}$ that $\Gamma_{1}\dot
{y}_{1}+\Gamma_{2}\dot{y}_{2}=0$.

For completeness we present the dynamical equations for the $x$-components,
for reference as needed when used.%

\begin{multline}
\dot{x}_{1}=-\left.  \frac{\partial\psi_{1}}{\partial y}\right\vert
_{_{\left(  x_{2},y_{2}\right)  }}=\Gamma_{1}\left[  K_{1}\left(
2y_{1}\right)  +\frac{1}{2y_{1}}\right]  \ +\frac{\Gamma_{2}\left(
y_{1}-y_{2}\right)  }{r_{1,2}}\left[  K_{1}\left(  r_{1,2}\right)  -\frac
{1}{r_{1,2}}\right] \label{dynamicx1}\\
+\frac{\Gamma_{2}\left(  y_{1}+y_{2}\right)  }{r_{1,2^{\ast}}}\left[  \frac
{1}{r_{1,2^{\ast}}}-K_{1}\left(  r_{1,2^{\ast}}\right)  \right]  \text{.}%
\end{multline}

\bigskip%

\begin{multline}
\dot{x}_{2}=-\left.  \frac{\partial\psi_{2}}{\partial y}\right\vert
_{_{\left(  x_{2},y_{2}\right)  }}=\Gamma_{1}\left[  K_{1}\left(
2y_{2}\right)  +\frac{1}{2y_{2}}\ \right]  +\frac{\Gamma_{2}\left(
y_{2}-y_{1}\right)  }{r_{1,2}}\left[  K_{1}\left(  r_{1,2}\right)  -\frac
{1}{r_{1,2}}\right] \label{dynamicx2}\\
+\frac{\Gamma_{2}\left(  y_{1}+y_{2}\right)  }{r_{1,2^{\ast}}}\left[  \frac
{1}{r_{1,2^{\ast}}}-K_{1}\left(  r_{1,2^{\ast}}\right)  \right]  \text{.}%
\end{multline}

\section{Equations of motion for symmetric integrable 3 vortex configurations}

We provide the equations of motion for the two cases considered.

\subsection{The case $\Gamma_{1}=$ $\Gamma_{2}=1,\Gamma_{3}=-\alpha$ with
$\Gamma_{1}$ at $\left(  x,y\right)  $, $\Gamma_{2}$ at $\left(  -x,y\right)
$ (in layer 1) and $\Gamma_{3}$ at $\left(  0,y_{3}\right)  $ (in layer 2)}

In this case it can be shown that $\dot{x}_{1}=\dot{x}_{2}$ and $\dot{y}%
_{1}=-\dot{y}_{2}$. Relative equilibrium solutions will then be admitted, for
this configuration by requiring that $\dot{y}_{1}=0$, and $\dot{x}_{1}=\dot
{x}_{3}$. In this case the invariant $\Gamma_{1}y_{1}+\Gamma_{2}y_{2}%
+\Gamma_{3}y_{3}=\beta$, becomes $y_{1}+y_{2}-\alpha y_{3}=\beta$, or
$y_{3}=\frac{y_{1}+y_{2}-\beta}{\alpha}$.%

\begin{align}
\dot{y}_{1}  &  =\frac{1\text{\ }}{2x}+K_{1}\left(  2x\right)  -\frac
{2x}{\left(  2x\right)  ^{2}+\left(  2y\right)  ^{2}}-\frac{2xK_{1}\left(
\sqrt{\left(  2x\right)  ^{2}+\left(  2y\right)  ^{2}}\right)  }{\sqrt{\left(
2x\right)  ^{2}+\left(  2y\right)  ^{2}}}\nonumber\\
&  -\frac{\alpha x}{x^{2}+\left(  y-y_{3}\right)  ^{2}}+\frac{\alpha
xK_{1}\left(  \sqrt{x^{2}+\left(  y-y_{3}\right)  ^{2}}\right)  }{\sqrt
{x^{2}+\left(  y-y_{3}\right)  ^{2}}}\nonumber\\
&  +\frac{\alpha x}{x^{2}+\left(  y+y_{3}\right)  ^{2}}-\frac{\alpha
xK_{1}\left(  \sqrt{x^{2}+\left(  y+y_{3}\right)  ^{2}}\right)  }{\sqrt
{x^{2}+\left(  y+y_{3}\right)  ^{2}}}%
\end{align}

\bigskip%

\begin{align}
\dot{x}_{1}  &  =-\frac{1\text{\ }}{2y}-K_{1}\left(  2y\right)  -\frac
{2y}{\left(  2x\right)  ^{2}+\left(  2y\right)  ^{2}}-\frac{2yK_{1}\left(
\sqrt{\left(  2x\right)  ^{2}+\left(  2y\right)  ^{2}}\right)  }{\sqrt{\left(
2x\right)  ^{2}+\left(  2y\right)  ^{2}}}\nonumber\\
&  -\frac{\alpha\left(  y-y_{3}\right)  }{x^{2}+\left(  y-y_{3}\right)  ^{2}%
}+\frac{\alpha\left(  y-y_{3}\right)  K_{1}\left(  \sqrt{x^{2}+\left(
y-y_{3}\right)  ^{2}}\right)  }{\sqrt{x^{2}+\left(  y-y_{3}\right)  ^{2}}%
}\nonumber\\
&  +\frac{\alpha\left(  y+y_{3}\right)  }{x^{2}+\left(  y+y_{3}\right)  ^{2}%
}-\frac{\alpha\left(  y+y_{3}\right)  K_{1}\left(  \sqrt{x^{2}+\left(
y+y_{3}\right)  ^{2}}\right)  }{\sqrt{x^{2}+\left(  y+y_{3}\right)  ^{2}}}%
\end{align}

\bigskip%

\begin{align}
\dot{x}_{3}  &  =\frac{\alpha\text{\ }}{2y_{3}}+\alpha K_{1}\left(
2y_{3}\right)  +\frac{2\left(  y_{3}-y\right)  }{x^{2}+\left(  y_{3}-y\right)
^{2}}-\frac{2\left(  y_{3}-y\right)  K_{1}\left(  \sqrt{x^{2}+\left(
y_{3}-y\right)  ^{2}}\right)  }{\sqrt{x^{2}+\left(  y_{3}-y\right)  ^{2}}%
}\nonumber\\
&  -\frac{2\left(  y_{3}+y\right)  }{x^{2}+\left(  y_{3}+y\right)  ^{2}}%
-\frac{2\left(  y_{3}+y\right)  K_{1}\left(  \sqrt{x^{2}+\left(
y_{3}+y\right)  ^{2}}\right)  }{\sqrt{x^{2}+\left(  y_{3}+y\right)  ^{2}}}%
\end{align}

\subsection{The case $\Gamma_{1}=$ $\Gamma_{2}=1,\Gamma_{3}=-\alpha$ with
$\Gamma_{1}$ at $\left(  0,y_{1}\right)  $, $\Gamma_{2}$ at $\left(
0,y_{2}\right)  $ (in layer 1) and $\Gamma_{3}$ at $\left(  0,y_{3}\right)  $
(in layer 2)}

In this case it can be shown that $\dot{y}_{1}=\dot{y}_{2}=\dot{y}_{3}=0$.
Relative equilibrium solutions will then be admitted, for this configuration
by requiring that $\dot{x}_{1}=\dot{x}_{2}=\dot{x}_{3}$. Again the invariant
$\Gamma_{1}y_{1}+\Gamma_{2}y_{2}+\Gamma_{3}y_{3}=\beta$, becomes $y_{1}%
+y_{2}-\alpha y_{3}=\beta$, or $y_{3}=\frac{y_{1}+y_{2}-\beta}{\alpha}$.

\bigskip%

\begin{align}
\dot{x}_{1}  &  =-\frac{1\text{\ }}{2y_{1}}-K_{1}\left(  2y_{1}\right)
+\frac{\left(  y_{1}-y_{2}\right)  }{\left\vert y_{1}-y_{2}\right\vert ^{2}%
}+\frac{\left(  y_{1}-y_{2}\right)  K_{1}\left(  \left\vert y_{1}%
-y_{2}\right\vert \right)  }{\left\vert y_{1}-y_{2}\right\vert }\nonumber\\
&  -\frac{1}{y_{1}+y_{2}}-K_{1}\left(  y_{1}+y_{2}\right) \nonumber\\
-  &  \frac{\alpha\left(  y_{1}-y_{3}\right)  }{\left\vert y_{1}%
-y_{3}\right\vert ^{2}}+\frac{\alpha\left(  y_{1}-y_{3}\right)  K_{1}\left(
\left\vert y_{1}-y_{3}\right\vert \right)  }{\left\vert y_{1}-y_{3}\right\vert
}+\frac{\alpha}{y_{1}+y_{3}}-\alpha K_{1}\left(  y_{1}+y_{3}\right)
\end{align}

\begin{align}
\dot{x}_{2}  &  =-\frac{1\text{\ }}{2y_{2}}-K_{1}\left(  2y_{2}\right)
+\frac{\left(  y_{2}-y_{1}\right)  }{\left\vert y_{2}-y_{1}\right\vert ^{2}%
}+\frac{\left(  y_{2}-y_{1}\right)  K_{1}\left(  \left\vert y_{2}%
-y_{1}\right\vert \right)  }{\left\vert y_{2}-y_{1}\right\vert }\nonumber\\
&  -\frac{1}{y_{1}+y_{2}}-K_{1}\left(  y_{1}+y_{2}\right) \nonumber\\
-  &  \frac{\alpha\left(  y_{2}-y_{3}\right)  }{\left\vert y_{2}%
-y_{3}\right\vert ^{2}}+\frac{\alpha\left(  y_{2}-y_{3}\right)  K_{1}\left(
\left\vert y_{2}-y_{3}\right\vert \right)  }{\left\vert y_{2}-y_{3}\right\vert
}+\frac{\alpha}{y_{2}+y_{3}}-\alpha K_{1}\left(  y_{2}+y_{3}\right)
\end{align}

\begin{align}
\dot{x}_{3}  &  =\frac{\alpha\text{\ }}{2y_{3}}+\alpha K_{1}\left(
2y_{3}\right) \nonumber\\
&  +\frac{\left(  y_{3}-y_{1}\right)  }{\left\vert y_{3}-y_{1}\right\vert
^{2}}-\frac{\left(  y_{3}-y_{1}\right)  K_{1}\left(  \left\vert y_{3}%
-y_{1}\right\vert \right)  }{\left\vert y_{3}-y_{1}\right\vert }-\frac
{1}{y_{1}+y_{3}}+K_{1}\left(  y_{1}+y_{3}\right) \nonumber\\
&  +\frac{\left(  y_{3}-y_{2}\right)  }{\left\vert y_{3}-y_{2}\right\vert
^{2}}-\frac{\left(  y_{3}-y_{2}\right)  K_{1}\left(  \left\vert y_{3}%
-y_{2}\right\vert \right)  }{\left\vert y_{3}-y_{2}\right\vert }-\frac
{1}{y_{2}+y_{3}}+K_{1}\left(  y_{2}+y_{3}\right)
\end{align}

\bigskip
\end{document}